\documentclass[preprint,12pt]{elsarticle}
\usepackage{graphics}
\usepackage{graphicx}
\usepackage{subfigure}

\usepackage{latexsym}
\usepackage{amsmath}
\usepackage{graphicx}
\usepackage{times}
\usepackage{graphicx,multicol,enumerate,subfigure,amsmath}

\usepackage{fullpage}
\usepackage{setspace}
\usepackage[english]{babel}
\usepackage[version=3]{mhchem}
\usepackage{amsmath, amssymb, amsthm}
\usepackage{verbatim, latexsym}
\usepackage{colortbl}
\usepackage{multirow}

\newtheorem{theorem}{Theorem}
\newtheorem{remark}[theorem]{Remark}

\usepackage{soul}
\setstcolor{red}

\biboptions{sort&compress}
\journal{Journal of Our Choice}

\onehalfspace

\begin{document}

\begin{frontmatter}

\title{Generalized Multiscale
Finite Element Methods. Oversampling Strategies}
\author{\textbf{Yalchin Efendiev}$^{1,2*}$}

\author{\textbf{Juan Galvis$^{2,3}$, Guanglian Li$^{2}$ and Michael Presho$^{2}$}}

\address{$^{1}$ Center for Numerical Porous Media (NumPor) \\
King Abdullah University of Science and Technology (KAUST) \\
Thuwal 23955-6900, Kingdom of Saudi Arabia.}

\address{$^{2}$ Department of Mathematics \& Institute for Scientific Computation (ISC) \\
Texas A\&M University \\
College Station, Texas, USA}

\address{$^{3}$ Departamento de Matem\'{a}ticas\\ Universidad Nacional de Colombia\\ 
Bogot\'a D.C., Colombia \\
}

\cortext[cor1]{Email address: efendiev@math.tamu.edu}

\begin{abstract}
In this paper, we propose oversampling strategies in the Generalized 
Multiscale Finite Element Method (GMsFEM) framework. 
The GMsFEM, which has been recently introduced in \cite{egh12}, allows solving
multiscale parameter-dependent problems at a reduced computational cost
by constructing a reduced-order representation of the solution on a coarse
grid. The main idea of the method consists of (1) the construction
of snapshot space, (2) the construction of the offline space, and (3)
construction of the online space (the latter for parameter-dependent problems).
 In \cite{egh12}, it was shown that the GMsFEM
provides a flexible tool to solve multiscale problems 
with a complex input space by generating
appropriate snapshot, offline, and online spaces. In this paper,
we develop oversampling techniques to be used in this context (see \cite{hw97} 
where
oversampling is introduced for multiscale finite element methods).
It is known (see \cite{hw97})
that the oversampling can improve the accuracy of
multiscale methods. In particular, the oversampling technique uses
larger regions (larger than the target coarse block) in constructing
local basis functions.
Our motivation stems from the analysis presented in this paper which show
that when using oversampling techniques in the construction of the snapshot 
space
and offline space, GMsFEM will converge independent of
small scales and high-contrast under certain assumptions. 
We consider the use
of multiple eigenvalue problem to improve the convergence
and discuss their relation to single spectral problems that use oversampled
regions.
 The oversampling procedures
proposed in this paper differ from those in \cite{hw97}. In particular,
the oversampling domains are partially used in constructing local
spectral problems. We present numerical results and compare various
oversampling techniques in order to complement the proposed technique and analysis.

\end{abstract}


\end{frontmatter}

\section{Introduction}

Heterogeneous media with multiple scales and high-contrast
commonly occur in many applications, such as porous media and material
sciences. 
The development of reduced-order models describing complex processes
in such media is needed in such applications.
There are a variety of multiscale methods, e.g.
 \cite{aarnes04,apwy07,eh09,hw97,hughes98,jennylt03},
that efficiently capture multiscale behavior by constructing a reduced
representation of the solution space on a coarse grid.
While standard multiscale methods have proven effective for a
variety of applications (see, e.g.,
\cite{eghe05,eh09,ehg04,jennylt03}), in this paper we
consider a more recent framework, GMsFEM, in which the coarse spaces may be
systematically enriched to converge to the fine grid
solution. In particular, we develop oversampling techniques within GMsFEM
and show that these methods converge independent
of the small scales and high contrast under certain assumptions.

The Generalized Multiscale Finite Element Method (GMsFEM) is a flexible
framework that generalizes the Multiscale Finite Element Method (MsFEM)
by systematically enriching the coarse spaces and taking into account
small scale information and
complex input spaces. This approach, as in many
multiscale model reduction techniques, divides the computation into
two stages: offline and online. In the offline stage, 
a small dimensional space is constructed that can be efficiently
used in the online stage to construct multiscale basis functions.
These multiscale basis functions can be re-used for any input parameter
to solve the problem on a coarse grid. The main idea behind the construction
of offline and online spaces is the selection of local spectral problems
and the selection of the snapshot space.
In \cite{egh12}, we propose several general strategies. In this paper,
our focus is on the development of oversampling strategies.

Oversampling techniques have been developed in the context of
multiscale finite element methods \cite{hw97}
as well as  upscaling methods \cite{cdgw03}. These techniques
use the local solutions in larger domains to construct multiscale
basis functions in the context of MsFEM.
 We borrow that main concept in this paper. In particular,
we use the space of snapshots in the oversampled regions
by constructing a snapshot space spanned by harmonic functions or 
dominant eigenvectors of a local spectral problem formulated
in the oversampled domain. Furthermore, we use special
local spectral problems to determine the dominant modes
in the space of snapshots. This spectral problem is motivated by
the analysis and it
 uses a weighted mass matrix in the oversampled region while
the energy (stiffness) matrix is constructed in the target coarse domain.
By choosing the dominant modes, we identify multiscale basis 
functions. These basis functions are then multiplied by partition
of unity functions to solve the flow equation on a coarse grid
(in the absence of the parameter). 
We also describe the use of multiple local spectral problems
for enhancing the accuracy of the approximation
and discuss their relation to single spectral problems that use oversampled
regions where the latter provides an optimal space.
In the presence of the parameter,
we also design an online space following the same strategy as the offline
space construction (but using an online parameter value).
We employ the Galerkin finite element method, though discontinuous Galerkin
methods can also be used \cite{DGcoupling}.

We present numerical results that demonstrate the convergence of the
proposed methods. In our numerical experiments,
we test two different snapshot spaces that consist
of harmonic functions in the oversampled region and 
dominant eigenmodes of a local spectral problem in the oversampled
region. For the local spectral problems, we also consider
various choices by considering mass and energy matrices in the oversampled
regions. Our numerical results show that the proposed methods converge as
we increase the dimension of the space and this convergence is consistent
with our theoretical findings. We also test the use of
multiple spectral problems in constructing basis functions
as well as modifying the conductivity outside the target block to
improve the accuracy.

The paper is organized in the following way. In the next section,
Section \ref{prelim}, we present the problem setting and the definitions
of coarse grids. In Section \ref{locbasis}, we present the construction of local
basis functions. Section \ref{sec:numex} is devoted to the numerical results. 
In Section \ref{appendix:convergence} we present the analysis of the method and in 
Section \ref{sec:conclusions} we offer some concluding remarks.

\section{Preliminaries}
\label{prelim}

\begin{figure}\centering
 \subfigure[$\kappa_1(x)$]{\label{fig:permi}
    \includegraphics[width = 0.45\textwidth]{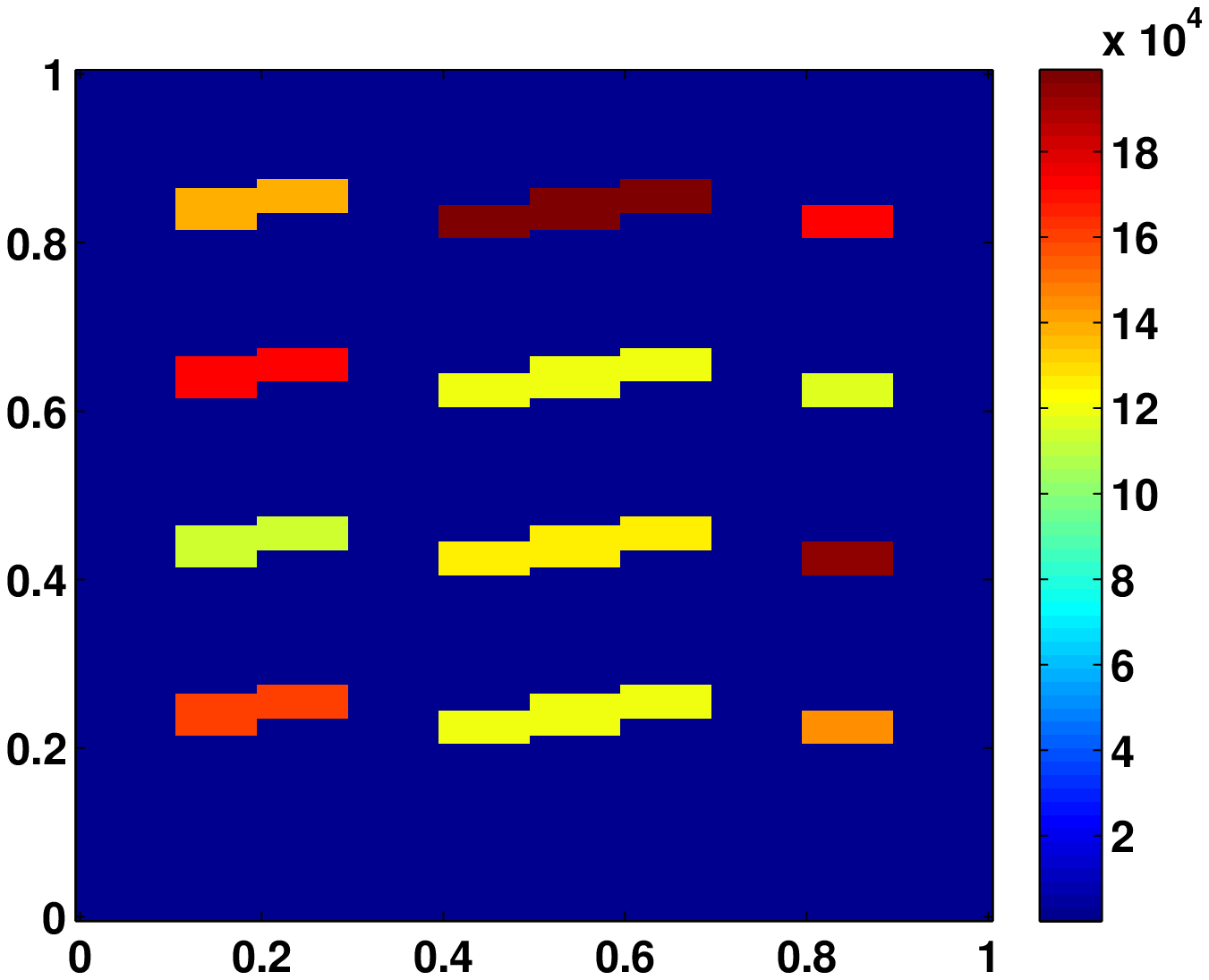}
   }
  \subfigure[$\kappa_2(x)$]{\label{fig:permii}
     \includegraphics[width = 0.45\textwidth]{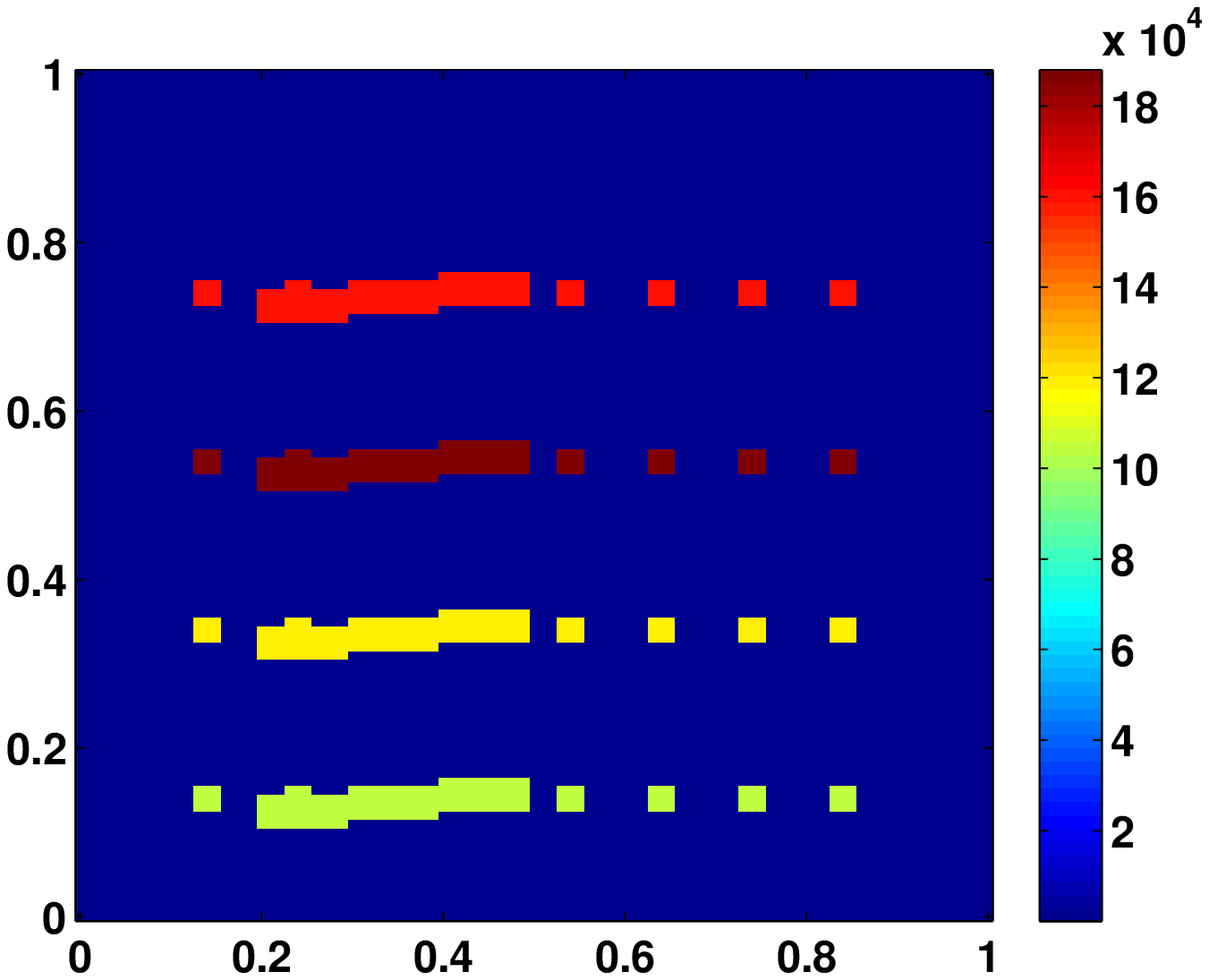}
  }
 \caption{Decomposition of permeability field}\label{perm}
\end{figure}

We consider elliptic equations of the form
\begin{equation} \label{eq:original}
-\mbox{div} \big( \kappa(x;\mu) \, \nabla u  \big)=f \, \, \text{in} \, D,
\end{equation}
where $u$ is prescribed on $\partial D$ and $\mu$ is a parameter. We assume
that $\kappa(x;\mu)=\sum_{q=1}^Q \Theta(\mu_q) \kappa_q(x)$ and
that the coefficient
$\kappa(x;\cdot)$ has multiple scale and high variations
(e.g., see Fig.~\ref{perm} for  $\kappa_1(x)$ and 
$\kappa_2(x)$ used in simulations).

\begin{figure}[tb]
  \centering
  \includegraphics[width=1.0 \textwidth]{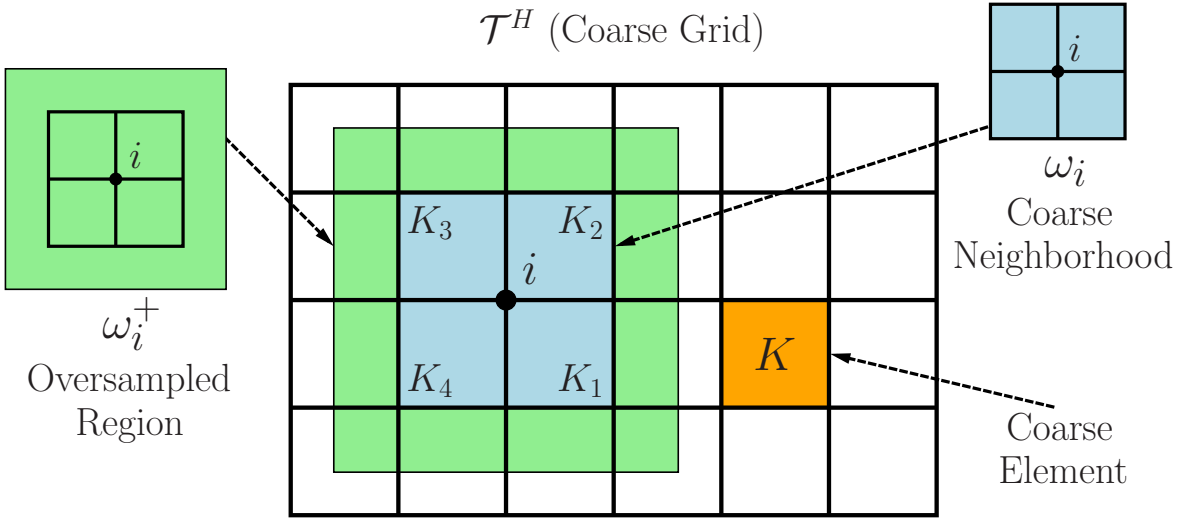}
  \caption{Illustration of a coarse neighborhood and oversampled domain}
  \label{schematic}
\end{figure}

To discretize (\ref{eq:original}), we next introduce the notion of fine and coarse grids. We let $\mathcal{T}^H$ be a usual conforming partition of the computational domain $D$ into finite elements (triangles, quadrilaterals, tetrahedrals, etc.). We refer to this partition as the coarse grid and assume that each coarse subregion is partitioned into a connected union of fine grid blocks. The fine grid partition will be denoted by $\mathcal{T}^h$. We use $\{x_i\}_{i=1}^{N_v}$ (where $N_v$ the number of coarse nodes) to denote the vertices of
the coarse mesh $\mathcal{T}^H$, and define the neighborhood of the node $x_i$ by
\begin{equation} \label{neighborhood}
\omega_i=\bigcup\{ K_j\in\mathcal{T}^H; ~~~ x_i\in \overline{K}_j\}.
\end{equation}
See Fig.~\ref{schematic} for an illustration of neighborhoods and elements subordinated to the coarse discretization. Furthermore, we introduce a notation for an oversampled region. We
denote by $\omega_i^{+}$ an oversampled region of $\omega_i\subset \omega_i^{+}$.
In general, we will consider oversampled regions $\omega_i^{+}$
defined by adding several fine-grid   or coarse-grid layers around $\omega_i$.

Next, we briefly outline the global coupling and the role of coarse basis
functions for the respective formulations that we consider. Throughout this paper, we use the continuous Galerkin formulation, and use $\omega_i$ as the support of basis functions even though $\omega_i^{+}$ will be used
in constructing multiscale basis functions.
 For the purpose of this description, we formally denote the basis functions of the online space $V_{\text{on}}$ by $\psi_k^{\omega_i}$. The solution will be sought as $u_H(x;\mu)=\sum_{i,k} c_{k}^i \psi_{k}^{\omega_i}(x; \mu)$.

Once the basis functions are identified, the global coupling is given through the variational form
\begin{equation}
\label{eq:globalG} a(u_H,v;\mu)=(f,v), \quad \text{for all} \, \, v\in
V_{\text{on}},
\end{equation}
and 
\[
a(u,v;\mu)=\int_D \kappa(x;\mu)\nabla u \nabla v.
\]
We note that in the case when the coefficient is independent of the parameter,
then $V_{\text{on}}=V_{\text{off}}$.

\section{Local basis functions}
\label{locbasis}
In this section we describe the offline-online computational procedure, and elaborate on some applicable choices for the associated bilinear forms to be used in the coarse space construction. Below we offer a general outline for the procedure.

\begin{itemize}
\item[1.]  Offline computations:
\begin{itemize}
\item 1.0. Coarse grid generation.
\item 1.1. Construction of snapshot space that will be used to compute an offline space.
\item 1.2. Construction of a small dimensional offline space by performing dimension reduction in the space of local snapshots.
\end{itemize}
\item[2.] Online computations:
\begin{itemize}
\item 2.1. For each input parameter, compute multiscale basis functions.
\item 2.2. Solution of a coarse-grid problem for any force term and boundary condition.
\item 2.3. Iterative solvers, if needed.
\end{itemize}
\end{itemize}

In the offline computation, we first construct a snapshot space $V_{\text{snap}}^{\omega_i^+}$ or $V_{\text{snap}}^{\omega_i}$, depending on the choice of domain to generate the snapshot space, where $\omega_i^+$ is an oversampled region that contains a coarse neighborhood $\omega_i$. Construction of the snapshot space involves solving the local problems for various choices of input parameters, and we describe the details below.

\subsection{Snapshot space}

\subsubsection{Harmonic extensions in oversampled region}
\label{sec:harmonic}
Our first choice of snapshot space consists of harmonic extension
of fine-grid functions defined on the boundary of $\omega_i^{+}$.
More precisely, for each fine-grid function, $\delta_l^h(x)$,
which is defined by 
$\delta_l^h(x)=\delta_{l,k},\,\forall l,k\in \textsl{J}_{h}(\omega_i^{+})$, where $\textsl{J}_{h}(\omega_i^{+})$ denotes the fine-grid boundary node on $\partial\omega_i^{+}$. 

For parameter-independent problem, we
solve
\[
-div(\kappa (x)\nabla  \psi_{l}^{+, \text{snap}})=0\ \ \text{in} \ \omega_i^{+}
\]
subject to boundary condition, $ \psi_{l}^{+, \text{snap}}=\delta_l^h(x)$
on $\partial\omega_i^{+}$. 

For parameter-dependent one, we can choose several values $\mu_j,\,j=1,\ldots,J$ ($J$ denotes the number of parameters used) to generate the snapshot space separately as above and combine them to obtain the snapshot space
(see details in Section \ref{sec:pdcase}).

\subsubsection{Local spectral basis in oversampled region}
\label{sec:spectral}
We propose to solve the following zero Neumann eigenvalue problem on an oversampled domain $\omega_i^+$:
\begin{equation} \label{snaplinalg}
A^+(\mu_j) \psi_{l,j}^{+, \text{snap}} = \lambda_{l,j}^{+, \text{snap}} S^+(\mu_j) \psi_{l,j}^{+, \text{snap}}
\,\, \text{in} \, \, \, \omega_i^+,
\end{equation}
where $\mu_j$ ($j=1,\dots,J$) is a specified set of fixed parameter values, and we emphasize that the superscript $+$ signifies that the eigenvalue problem is solved in an oversampled coarse subdomain $\omega_i^+$. The matrices in Eq.~\eqref{snaplinalg} are defined as
\begin{equation}\label{eqn:eigenmatrix}
A^+(\mu_j) = [a^+(\mu_j)_{mn}] = \int_{\omega_i^+} \kappa(x; \mu_j) \nabla \phi_n \cdot \nabla \phi_m \quad \text{and}
\quad S^+(\mu_j) = [s^+(\mu_j)_{mn}] = \int_{\omega_i^+} \widetilde{\kappa}(x; \mu_j)  \phi_n \phi_m,
\end{equation}
where $\phi_n$ denotes the standard bilinear, fine-scale basis
functions and the form for $\widetilde{\kappa}$ will be discussed in Section
\ref{appendix:convergence}. In our numerical implementations, we take 
$\widetilde{\kappa}=\kappa$, though one can use multiscale basis functions,
$\chi_i^{+}$
in $\omega_i^{+}$, to construct $\widetilde{\kappa}$ as
$\widetilde{\kappa} =\sum_i \kappa |\nabla \chi_i^{+}|^2$
(see \cite{ge09_3, egt11} for more discussions on the choice
of partition of unity functions).
 We note that Eq.~\eqref{snaplinalg} is the discretized
form of the continuous equation
\begin{equation*}
-\text{div}(\kappa(x; \mu_j) \nabla \psi_{l,j}^{+, \text{snap}} ) = \lambda_{l,j}^{+, \text{snap}} \widetilde{\kappa}(x; \mu_j)  \psi_{l,j}^{+, \text{snap}}
\quad \text{in} \, \, \, \omega_i^+.
\end{equation*}

After solving Eq.~\eqref{snaplinalg}, we keep the first $L_i$ eigenfunctions corresponding to the dominant eigenvalues (asymptotically vanishing in this case) to form the space
$$
V_{\text{snap}}^+ = \text{span}\{ \psi_{l,j}^{+, \text{snap}}:~~~1\leq j \leq J ~~ \text{and} ~~ 1\leq l \leq L_i \},
$$
for each oversampled coarse neighborhood $\omega_i^+$. We note that in the case when $\omega_i$ is adjacent to the global boundary, no oversampled domain is used. For the sake of simplicity, throughout,
 we denote continuous and discrete solutions
by the same symbol (e.g., $\psi_{l,j}^{+, \text{snap}}$ in the above case).

 We reorder the snapshot functions using a single index to create the matrices
$$
R_{\text{snap}}^+ = \left[ \psi_{1}^{+,\text{snap}}, \ldots, \psi_{M_{\text{snap}}}^{+,\text{snap}} \right] \quad \text{and} \quad R_{\text{snap}} = \left[ \psi_{1}^{\text{snap}}, \ldots, \psi_{M_{\text{snap}}}^{\text{snap}} \right],
$$
where $\psi_j^{\text{snap}}$ denotes the restriction of $\psi_j^{\text{+,snap}}$ to $\omega_i$, and $M_{\text{snap}}$ denotes the total number of functions to keep in the snapshot matrix construction.

Note that the above process to generate local spectral basis is also
 applied to parameter-independent problems. 

\subsection{Offline space}

We will discuss two types of offline spaces where the first one will
use one spectral problem in the snapshot space and the other one
will use multiple spectral problems in the snapshot space (following
Theorem 3.3 of \cite{bl11}).

\subsubsection{Offline space using a single spectral problem}

In order to construct an oversampled offline space $V_{\text{off}}^+$ or standard neighborhood offline space $V_{\text{off}}$, we perform a dimension reduction in the space of snapshots using an auxiliary spectral decomposition. The main objective is to use the offline space to efficiently (and accurately) construct a set of multiscale basis functions for each $\mu$ value in the online stage. More precisely, we seek a subspace of the snapshot space such that it can approximate any element of the snapshot space in the appropriate sense defined via auxiliary bilinear forms. At the offline stage the bilinear forms are chosen to be \emph{parameter-independent}, such that there is no need to reconstruct the offline space for each $\mu$ value. We will consider the following eigenvalue problems in the space of snapshots:

\begin{eqnarray}
A^{\text{off}} \Psi_k^{\text{off}} &=& \lambda_k^{\text{off}} S^{\text{off}}\Psi_k^{\text{off}}  \label{offeig1} \\
 &\text{or}   \nonumber \\
A^{+,\text{off}}\Psi_k^{\text{off}} &=& \lambda_k^{\text{off}} A^{\text{off}} \Psi_k^{\text{off}}  
\label{offeig2} \\ 
 &\text{or}  \nonumber \\
A^{\text{off}}\Psi_k^{\text{off}} &=& \lambda_k^{\text{off}} S^{+, \text{off}} \Psi_k^{\text{off}}
\label{offeig3} \\
 &\text{or}  \nonumber \\
A^{+,\text{off}}\Psi_k^{\text{off}} &=& \lambda_k^{\text{off}} S^{+, \text{off}} \Psi_k^{\text{off}}
\label{offeig4}
\end{eqnarray}
where
\begin{equation*}
 \displaystyle A^{\text{off}}= [a^{\text{off}}_{mn}] = \int_{\omega_i} \overline{\kappa}(x; \mu) \nabla \psi_m^{\text{snap}} \cdot \nabla \psi_n^{\text{snap}} = R_{\text{snap}}^T \overline{A} R_{\text{snap}},
 \end{equation*}
\begin{equation*}
 \displaystyle S^{\text{off}} = [s^{\text{off}}_{mn}] = \int_{\omega_i} \widetilde{\overline{\kappa}}(x; \mu) \psi_m^{\text{snap}} \psi_n^{\text{snap}} = R_{\text{snap}}^T \overline{S} R_{\text{snap}},
 \end{equation*}
\begin{equation*}
 \displaystyle A^{+,\text{off}} = [a_{mn}^{+,\text{off}}] = \int_{\omega_i^+} \overline{\kappa}(x,\mu) \nabla \psi_m^{+,\text{snap}} \cdot \nabla \psi_n^{+,\text{snap}} = \big(R_{\text{snap}}^+\big)^T \overline{A}^+ R_{\text{snap}}^+,
 \end{equation*}
 \begin{equation*}
 \displaystyle S^{+,\text{off}} = [s_{mn}^{+,\text{off}}] = \int_{\omega_i^+} \widetilde{\overline{\kappa}}(x,\mu) \psi_m^{+,\text{snap}}  \psi_n^{+,\text{snap}} = \big(R_{\text{snap}}^+\big)^T \overline{S}^+ R_{\text{snap}}^+.
 \end{equation*}
The coefficients $\overline{\kappa}(x, \mu) $ and $\widetilde{\overline{\kappa}}(x, \mu)$ are parameter-averaged coefficients (see \cite{egh12}).
Again, we will take  $\widetilde{\overline{\kappa}}(x, \mu)=\overline{\kappa}(x, \mu)$ though one can use multiscale partition of
unity functions to compute  $\widetilde{\overline{\kappa}}(x, \mu)$
(cf. \cite{egt11}).
 We note that $\overline{A}^+$ and $\overline{A}$ denote analogous fine scale matrices as defined in Eq.~\eqref{snaplinalg}, except that 
parameter-averaged coefficients are used in the construction, and that $A$ is constructed by integrating only on $\omega_i$. To generate the offline space we then choose the smallest $M_{\text{off}}$ eigenvalues from one of Eqs.~\eqref{offeig1}-\eqref{offeig3} and form the corresponding eigenvectors in the respective space of snapshots by setting
$\psi_k^{+,\text{off}} = \sum_j \Psi_{kj}^{\text{off}} \psi_j^{+,\text{snap}}$ or $\psi_k^{\text{off}} = \sum_j \Psi_{kj}^{\text{off}} \psi_j^{\text{snap}}$ (for $k=1,\ldots, M_{\text{off}}$), where $\Psi_{kj}^{\text{off}}$ are the coordinates of the vector $\Psi_{k}^{\text{off}}$. We then create the offline matrices
 $$
R_{\text{off}}^+ = \left[ \psi_{1}^{+,\text{off}}, \ldots, \psi_{M_{\text{off}}}^{+,\text{off}} \right] 
\quad \text{and} \quad R_{\text{off}} = \left[ \psi_{1}^{\text{off}}, \ldots, \psi_{M_{\text{off}}}^{\text{off}} \right] 
$$
to be used in the online space construction.

\begin{remark}
At this stage, we note that in the case when we have a parameter-independent coefficient in Eq.~\eqref{eq:original}, many of the expressions in this section are simplified. In particular, there is no need for averaging the coefficients in order to create the respective (offline) mass and  (offline) stiffness matrices. Furthermore, the offline space represents the final space in which the enriched multiscale solutions will be computed. Thus, the discussion of online space creation below is limited to the case when the problem is parameter-dependent. 
\end{remark}

\begin{figure}\centering
\subfigure[Eigenvalue problem Eq.~\eqref{offeig1} with Harmonic snapshots]{\label{fig:eig1}
     \includegraphics[width = 0.45\textwidth]{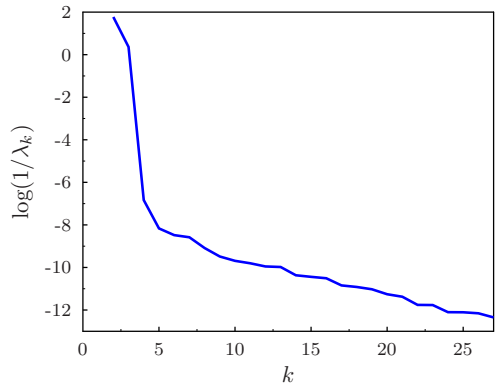}
  }
 \subfigure[Eigenvalue problem Eq.~\eqref{offeig2} with Harmonic snapshots]{\label{fig:eig2}
    \includegraphics[width = 0.45\textwidth]{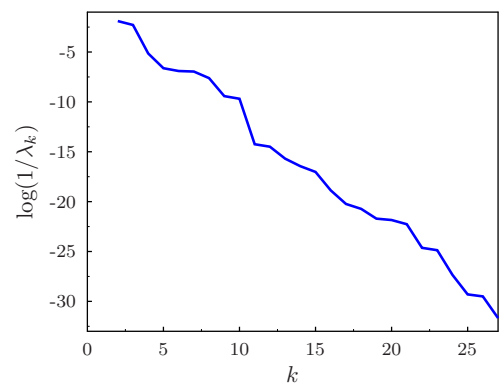}
   }
  \subfigure[Eigenvalue problem Eq.~\eqref{offeig3} with Harmonic snapshots]{\label{fig:eig3}
     \includegraphics[width = 0.45\textwidth]{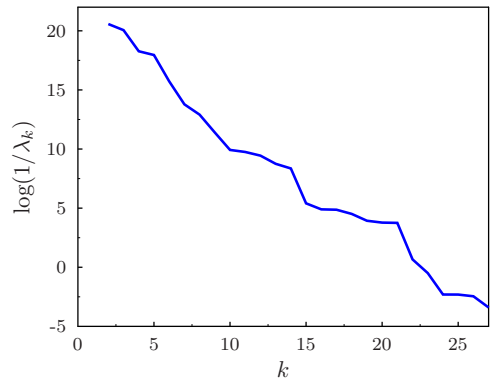}
  }
  \subfigure[Eigenvalue problem Eq.~\eqref{offeig3} with spectral snapshots]{\label{fig:eig4}
     \includegraphics[width = 0.45\textwidth]{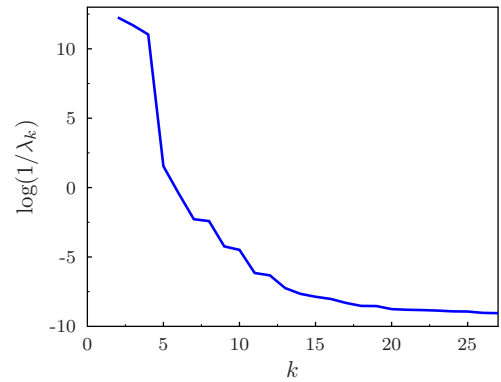}
  }
 \caption{Eigenvalue decay on log-scale against the number of eigenvalues. x-axis indicates the number of eigenvalue and y-axis indicates the inverse of the eigenvalue (on log-scale)}\label{fig:eigbehaviour}
\end{figure}

\begin{remark}
Our analysis in Section \ref{appendix:convergence} 
shows that the convergence of the GMsFEM is proportional to
the reciprocal of the eigenvalue that the corresponding eigenvector is not
included in the coarse space. We have compared the decay of the reciprocal
of eigenvalues for Eq.~\eqref{offeig1}, Eq.~\eqref{offeig2}, and Eq.~\eqref{offeig3} (by choosing a subdomain for $\kappa(x)$ in Fig.~\ref{fig:HCCoeff}).
We plot the decay of the eigenvalues for a coarse block in 
Fig.~\ref{fig:eigbehaviour} (note logarithmic y-scale).
As we observe from this figure that the decay of eigenvalues
corresponding to Eq.~\eqref{offeig3} (when oversampling is used in
formulating the eigenvalue problem) is faster compared to
 Eq.~\eqref{offeig1} (when no oversampling is used).
\end{remark}

\subsubsection{Offline space using multiple spectral problems}
\label{sec:multiple}

Motivated by Theorem 3.3 of \cite{bl11}, 
we propose an offline space that uses both Eq.~\eqref{offeig1}
and Eq.~\eqref{offeig4}. In particular, using dominant eigenvectors
of both Eq.~\eqref{offeig1}
and Eq.~\eqref{offeig4}, we take a union of these eigenvectors
to construct an offline space. In particular, as described above,
we use
$\psi_k^{+,\text{off}} = \sum_j \Psi_{kj}^{+,\text{off}} \psi_j^{+,\text{snap}}$ 
 (for $k=1,\ldots, M_{+,\text{off}}$)
or $\psi_k^{\text{off}} = \sum_j \Psi_{kj}^{\text{off}} \psi_j^{\text{snap}}$ (for $k=1,\ldots, M_{\text{off}}$), where $\Psi_{kj}^{+,\text{off}}$ are the coordinates of the vector $\Psi_{k}^{\text{off}}$ in Eq.~\eqref{offeig4} and  $\Psi_{kj}^{\text{off}}$ are the coordinates of the vector $\Psi_{k}^{\text{off}}$ in Eq.~\eqref{offeig1}. Then, the offline space is constructed as a union of $\psi_k^{\text{off}}$
and $\psi_k^{+,\text{off}}$ after eliminating linearly dependent vectors.
We present an analysis in Section \ref{appendix:convergence_multiple} and
 numerical results in Section \ref{sec:num_res_multiple}.

\subsection{Online space for parameter-dependent case}

We only describe the online space using a single spectral problem.
One can analogously construct the online space using multiple 
spectral problems.
For the parameter-dependent case, we next construct the associated online coarse space
$V_{\text{on}}(\mu)$ \emph{for each fixed} $\mu$ value on each coarse subdomain.
In principle, we want this to be a small dimensional subspace of the offline space for computational efficiency.
The online coarse space will be used within the finite element
framework to solve the original global problem, where a continuous Galerkin coupling of the multiscale basis functions is used to compute the global solution. In particular, we seek a subspace of the respective offline space such that it can approximate any element of the offline space in
an appropriate sense. We note that at the online stage, the bilinear forms are chosen to be \emph{parameter-dependent}. Similar analysis motivates the following eigenvalue problems posed in the offline space:
\begin{eqnarray}
A^{\text{on}}(\mu) \Psi_k^{\text{on}} &=& \lambda_k^{\text{on}} S^{\text{on}}(\mu)\Psi_k^{\text{on}}  \label{oneig1} \\
 &\text{or}   \nonumber \\
A^{+,\text{on}}(\mu)\Psi_k^{\text{on}} &=& \lambda_k^{\text{on}} A^{\text{on}} (\mu)\Psi_k^{\text{on}}  
\label{oneig2} \\ 
 &\text{or}  \nonumber \\
A^{\text{on}}(\mu)\Psi_k^{\text{on}} &=& \lambda_k^{\text{on}} S^{+, \text{on}}(\mu) \Psi_k^{\text{on}}
\label{oneig3}
\end{eqnarray}
where
\begin{equation*}
 \displaystyle A^{\text{on}}(\mu) = [a^{\text{on}}(\mu)_{mn}] = \int_{\omega_i} \kappa(x; \mu) \nabla \psi_m^{\text{off}} \cdot \nabla \psi_n^{\text{off}} = R_{\text{off}}^T A(\mu) R_{\text{off}},
 \end{equation*}
\begin{equation*}
 \displaystyle S^{\text{on}}(\mu) = [s^{\text{on}}(\mu)_{mn}] = \int_{\omega_i} \widetilde{\kappa}(x; \mu) \psi_m^{\text{off}} \psi_n^{\text{off}} = R_{\text{off}}^T S(\mu) R_{\text{off}},
 \end{equation*}
\begin{equation*}
 \displaystyle A^{+,\text{on}}(\mu) = [a_{mn}^{+,\text{on}}(\mu)] = \int_{\omega_i^+} \kappa(x,\mu) \nabla \psi_m^{+,\text{off}} \cdot \nabla \psi_n^{+,\text{off}} = \big(R_{\text{off}}^+\big)^T A^+(\mu) R_{\text{off}}^+,
 \end{equation*}
 \begin{equation*}
 \displaystyle S^{+,\text{on}}(\mu) = [s_{mn}^{+,\text{on}}(\mu)] = \int_{\omega_i^+} \widetilde{\kappa}(x,\mu) \psi_m^{+,\text{off}} \psi_n^{+,\text{off}} = \big(R_{\text{off}}^+\big)^T S^+(\mu) R_{\text{off}}^+,
 \end{equation*}
and $\kappa(x; \mu)$ and $\widetilde{\kappa}(x; \mu)$ are now parameter dependent. 
Again, we will take  $\widetilde{\kappa}(x, \mu)={\kappa}(x, \mu)$ in our simulations though one can use multiscale partition of
unity functions to compute  $\widetilde{\kappa}(x, \mu)$
(cf. \cite{egt11}).
To generate the online space we then choose the smallest $M_{\text{on}}$ eigenvalues from one of Eqs.~\eqref{oneig1}-\eqref{oneig3} and form the corresponding eigenvectors in the offline space by setting
$\psi_k^{\text{on}} = \sum_j \Psi_{kj}^{\text{on}} \psi_j^{\text{off}}$ (for $k=1,\ldots, M_{\text{on}}$), where $\Psi_{kj}^{\text{on}}$ are the coordinates of the vector $\Psi_{k}^{\text{on}}$.

\section{Numerical Examples}
\label{sec:numex}
\subsection{Parameter-independent case}\label{sec:parameter-independent}

\begin{figure}\centering
 
     \includegraphics[width = 0.50\textwidth, keepaspectratio = true]{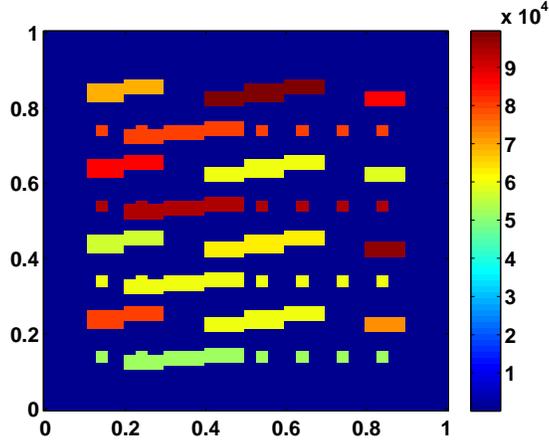}
  
 \caption{Permeability field used in Section \ref{sec:parameter-independent}}\label{fig:HCCoeff}
\end{figure}
First, we consider parameter-independent case
\[
\kappa(x;\mu)=\kappa(x)
\]
 by choosing
$\mu_1, \mu_2=0.5$ (see Fig. \ref{fig:HCCoeff} for an illustration of the resulting permeability). In previous works, e.g., \cite{egt11}, 
$\widetilde{\kappa}$ takes the general form $\widetilde{\kappa} = \kappa \sum_i H^2 | \nabla \chi_i |^2$, where $\chi_i$ denotes an original partition of unity \cite{egh12}, although we take $\widetilde{\kappa} = \kappa$ for the majority of examples in this section. The fine-grid is chosen to be
$100 \times 100$. We consider two coarse grids, $10\times 10$
and $20 \times 20$. The error will be measured in weighted $L^2$
and weighted $H^1$ norms defined as
\[
\|u\|_{L^2_\kappa}=\left(\int_D \kappa u^2\right)^{1\over 2},
\ \ \|u\|_{H^1_\kappa}=\left(\int_D \kappa |\nabla u|^2\right)^{1\over 2}.
\]
Then Eq.~\eqref{eq:original} is solved with $f=1$ and linear Dirichlet boundary
condition.

In the first set of numerical examples, 
$10\times 10$ coarse grid and the
oversampling region of the size of 10 fine-grid blocks
in each direction is chosen 
(i.e., the oversampled region contains an extra coarse block layer
around $\omega_i$). We denote this oversampled
region by $\omega_i^{+}=3\times \omega_i$.
We use the eigenvalue problems Eq.~\eqref{offeig1}, Eq.~\eqref{offeig2},
and Eq.~\eqref{offeig3} in the space of snapshots generated in
the oversampled region by harmonic extensions as in Subsection
\ref{sec:harmonic}. 
In all numerical cases, we take $\widetilde{\kappa}=\kappa$.
In Tables \ref{off10_parindep_At_St}, \ref{off10_parindep_Ap_At}, and \ref{off10_parindep_At_Sp},
we present the errors for weighted $L^2$ norm and weighted $H^1$ norm.
As we observe that all cases predict similar convergence errors
that decrease as we increase the dimension of the space. We note
that in this case there is a residual error because of the fact that
we use harmonic functions as a space of snapshots and thus
we can not approximate the error due to the source term. This error
for $10\times 10$ coarse mesh is about $10 \%$ (or order of coarse mesh size).
 Because of this irreducible
error, the convergence of GMsFEM deteriorates and remains at $10 \%$.

\begin{table}[htb]
\centering

\begin{tabular}{|c|c|c|c|c|c|}
\hline 
\multirow{2}{*}{$\text{dim}(V_{\text{off}})$} &\multirow{2}{*}{$\Lambda_*$} &
\multicolumn{2}{c|}{  $\|u-u^{\text{off}} \|$ (\%) }  \\
\cline{3-4} {}&{}&
$\hspace*{0.8cm}   L^{2}_\kappa(D)   \hspace*{0.8cm}$ &
$\hspace*{0.8cm}   H^{1}_\kappa(D)  \hspace*{0.8cm}$
\\
\hline\hline
       $364$     &     $1.04{\times}10^3$   &  $0.46$    & $17.26$  \\
\hline
      $526$    &   $2.62{\times}10^3$  & $0.42$    & $13.84$ \\
\hline
      $688$    &   $4.29{\times}10^3$  & $0.33$    & $11.92$ \\
\hline
       $909$   &  $1.15{\times}10^5$ & $0.30$  &$10.90$  \\

\hline
\end{tabular}

\caption{Relative errors between the fine scale solution and {\bf offline spaces}; Eigenvalue problem from Eq.~\eqref{offeig1}, $10{\times}10$ coarse mesh, harmonic snapshots, $\omega_i^{+}=3\times \omega_i$.}
\label{off10_parindep_At_St}
\end{table}

\begin{table}[htb]
\centering

\begin{tabular}{|c|c|c|c|c|c|}
\hline 
\multirow{2}{*}{$\text{dim}(V_{\text{off}}^+)$} &\multirow{2}{*}{$\Lambda_*$} &
\multicolumn{2}{c|}{  $\|u-u^{+, \text{off}} \|$ (\%) }  \\
\cline{3-4} {}&{}&
$\hspace*{0.8cm}   L^{2}_\kappa(D)   \hspace*{0.8cm}$ &
$\hspace*{0.8cm}   H^{1}_\kappa(D)  \hspace*{0.8cm}$
\\
\hline\hline
       $364$     &     $33.16$   &  $4.22$    & $33.69$  \\
\hline
      $526$    &   $105.32$  & $0.65$    & $16.59$ \\
\hline
      $688$    &   $669.50$  & $0.32$    & $11.95$ \\
\hline
       $909$   &  $8.12{\times}10^6$ & $0.30$  &$10.82$  \\

\hline
\end{tabular}

\caption{Relative errors between the fine scale solution and {\bf offline spaces}; Eigenvalue problem from
              Eq.~\eqref{offeig2}, $10{\times}10$ coarse mesh, harmonic snapshots, $\omega_i^{+}=3\times \omega_i$.}
\label{off10_parindep_Ap_At}
\end{table}

\begin{table}[htb]
\centering

\begin{tabular}{|c|c|c|c|c|c|}
\hline 
\multirow{2}{*}{$\text{dim}(V_{\text{off}}^+)$} &\multirow{2}{*}{$\Lambda_*$} &
\multicolumn{2}{c|}{  $\|u-u^{+, \text{off}} \|$ (\%) }  \\
\cline{3-4} {}&{}&
$\hspace*{0.8cm}   L^{2}_\kappa(D)   \hspace*{0.8cm}$ &
$\hspace*{0.8cm}   H^{1}_\kappa(D)  \hspace*{0.8cm}$
\\
\hline\hline
       $364$     &     $2.11{\times}10^{-4}$   &  $4.47$    & $45.99$  \\
\hline
      $526$    &   $8.62{\times}10^{-4}$  & $1.62$    & $27.65$ \\
\hline
      $688$    &   $0.0018$  & $0.28$    & $14.54$ \\
\hline
       $909$   &  $0.0093$ & $0.30$  &$11.09$  \\

\hline
\end{tabular}

\caption{Relative errors between the fine scale solution and {\bf offline spaces}; Eigenvalue problem from
              Eq.~\eqref{offeig3}, $10{\times}10$ coarse mesh, harmonic snapshots, $\omega_i^{+}=3\times \omega_i$.}

\label{off10_parindep_At_Sp}
\end{table}

In the next example, we consider a smaller oversampled region that includes
only one fine grid block.   We denote this by $\omega_i^{+}=\omega_i + 1$.
We have tested various 
oversampled region sizes and include only one representative example.
In this example (see results in Table \ref{off10_parindep_At_Stt1+} 
),
 we observe similar error behavior as those in previous examples.

\begin{table}[htb]
\centering
\begin{tabular}{|c|c|c|c|c|c|}
\hline 
\multirow{2}{*}{$\text{dim}(V_{\text{off}}^{+})$} &\multirow{2}{*}{$\Lambda_*$} &
\multicolumn{2}{c|}{  $\|u-u^{+,\text{off}} \|$ (\%) }  \\
\cline{3-4} {}&{}&
$\hspace*{0.8cm}   L^{2}_\kappa(D)   \hspace*{0.8cm}$ &
$\hspace*{0.8cm}   H^{1}_\kappa(D)  \hspace*{0.8cm}$
\\
\hline\hline
       $526$     &     $0.03$   &  $1.64$    & $25.51$  \\
\hline
      $850$    &   $0.08$  & $0.38$    & $14.82$ \\
\hline
      $2470$    &   $533.50$  & $0.31$    & $12.10$ \\
\hline
       $3280$   &  $1.08{\times}10^3$ & $0.30$  &$11.42$  \\

\hline
\end{tabular}

\caption{Relative errors between the fine scale solution and {\bf offline spaces}; Eigenvalue problem from Eq.~\eqref{offeig3}, $10{\times}10$ coarse mesh, harmonic snapshots, $\omega_i^{+}=\omega_i + 1$.}
\label{off10_parindep_At_Stt1+}
\end{table}

As we discussed earlier, the error between the fine scale solution 
and GMsFEM solution contains an irreducible error because of the fact
that the harmonic snapshots are used and these snapshots
can not approximate the effects of the right hand side.
This error can be easily estimated for high-contrast problems
considered in this paper
and it is of order $O(H)$. First, we consider the use
of dominant eigenvectors as a space of snapshots
 that correspond to smallest eigenvalues 
of Eq.~\eqref{offeig1} in $\omega_i^{+}=3\times \omega_i$ as a snapshot 
space. In this snapshot space, we apply Eq.~\eqref{offeig3} and
identify dominant modes in the target domain as before. 
The numerical
results are presented in Table \ref{off10eigsnap_parindep_At_Sp}.
As we observe from these results that the error is smaller when
eigenvector snapshots are used. In general, when comparing to the
fine-scale solution, one can also use fewer modes 
corresponding to the space of harmonic snapshots and some extra modes
that represent source term within the local domain (e.g., modes that
correspond to homogeneous Dirichlet eigenvalue problem). 
In Table \ref{off10_parindep_At_Sp_vs909}, we
present numerical results, where the GMsFEM solution is compared to the
solution computed in the space of harmonic snapshots. 
In this setup, there is no irreducible error and the method converges
to the fine scale solution. Moreover, we notice that the errors
are smaller.

\begin{table}[htb]
\centering

\begin{tabular}{|c|c|c|c|c|c|}
\hline 
\multirow{2}{*}{$\text{dim}(V_{\text{off}}^+)$} &\multirow{2}{*}{$\Lambda_*$} &
\multicolumn{2}{c|}{  $\|u-u^{+, \text{off}} \|$ (\%) }  \\
\cline{3-4} {}&{}&
$\hspace*{0.8cm}   L^{2}_\kappa(D)   \hspace*{0.8cm}$ &
$\hspace*{0.8cm}   H^{1}_\kappa(D)  \hspace*{0.8cm}$
\\
\hline\hline
       $364$     &     $0.0045$   &  $0.27$    & $17.49$  \\
\hline
      $688$    &   $0.055$  & $0.08$    & $9.88$ \\
\hline
      $1012$    &   $0.91$  & $0.07$    & $7.33$ \\
\hline
       $1660$   &  $37.3$ & $0.03$  &$4.04$  \\
\hline
       $3280$   &  $1.68{\times}10^3$ & $0.004$  &$1.10$  \\

\hline
\end{tabular}

\caption{Relative errors between the fine scale solution and {\bf offline spaces}; Eigenvalue problem from
              Eq.~\eqref{offeig3}, $10{\times}10$ coarse mesh, eigenvalue snapshots, $\omega_i^+=3\times \omega_i$.}

\label{off10eigsnap_parindep_At_Sp}
\end{table}

\begin{table}[htb]
\centering

\begin{tabular}{|c|c|c|c|c|c|}
\hline 
\multirow{2}{*}{$\text{dim}(V_{\text{off}}^+)$} &\multirow{2}{*}{$\Lambda_*$} &
\multicolumn{2}{c|}{  $\|u^{+,909}-u^{+, \text{off}} \|$ (\%) }  \\
\cline{3-4} {}&{}&
$\hspace*{0.8cm}   L^{2}_\kappa(D)   \hspace*{0.8cm}$ &
$\hspace*{0.8cm}   H^{1}_\kappa(D)  \hspace*{0.8cm}$
\\
\hline\hline
       $364$     &     $2.11{\times}10^{-4}$   &  $4.69$    & $44.36$  \\
\hline
      $526$    &   $8.62{\times}10^{-4}$  & $1.75$    & $25.17$ \\
\hline
      $688$    &   $0.0018$  & $0.16$    & $9.34$ \\
\hline
       $860$   &  $0.0070$ & $0.05$  &$3.85$  \\

\hline
\end{tabular}

\caption{Relative errors between the maximal dimension offline solution ($u^{+,909}$) and {\bf offline spaces}    
              obtained from 
              using oversampled domains; Eigenvalue problem from
              Eq.~\eqref{offeig3}, $10{\times}10$ coarse mesh, harmonic snapshots, $\omega_i^{+}=3\times \omega_i$.}

\label{off10_parindep_At_Sp_vs909}
\end{table}

For the next set of numerical examples, we use $20 \times 20 $ coarse-grid.
In Table \ref{off20_parindep_At_Sp}, we present numerical
results when the eigenvalue problem  Eq.~\eqref{offeig3} is used. 
As in the case of $10 \times 10$ coarse grid, there is an irreducible error;
however, it is lower (about $5\%$), because of the coarse mesh size. 
To remove the irreducible error, we compare the GMsFEM solution
to the solution computed with snapshot vectors in Table \ref{off20_parindep_At_Sp_vs3471}. As we observe that the error is smaller and it will converge to
zero as we increase the dimension of the coarse space. We also present an
error when a different oversampling domain size is used in Table 
\ref{off20_parindep_At_Stt1+}. The results are not sensitive to the
oversampling domain size as these results show.
In Table 
\ref{off20_parindep_At_Sp_eigen_t5}, we present relative errors
when the snapshot space is chosen to consist of eigenvectors
as defined in Eq.~\ref{snaplinalg} (cf. Table \ref{off10eigsnap_parindep_At_Sp}). 
In this case, similar to Table \ref{off10eigsnap_parindep_At_Sp}, we observe smaller
errors when the snapshot space consists of eigenvectors 
in Eq.~\ref{snaplinalg}. We also present a numerical
result in Table \ref{off20_parindep_At_SpChoppedCoeff}
where the coefficients in $\omega_i^{+}\backslash\omega_i$ reduced
by $1e+4$ to diminish the constant in the estimates presented in Section
\ref{appendix:convergence_single}.

\begin{table}[htb]
\centering

\begin{tabular}{|c|c|c|c|c|c|}
\hline 
\multirow{2}{*}{$\text{dim}(V_{\text{off}}^+)$} &\multirow{2}{*}{$\Lambda_*$} &
\multicolumn{2}{c|}{  $\|u-u^{+, \text{off}} \|$ (\%) }  \\
\cline{3-4} {}&{}&
$\hspace*{0.8cm}   L^{2}_\kappa(D)   \hspace*{0.8cm}$ &
$\hspace*{0.8cm}   H^{1}_\kappa(D)  \hspace*{0.8cm}$
\\
\hline\hline
       $1524$     &     $3.25{\times}10^{-5}$   &  $2.69$    & $36.37$  \\
\hline
      $2168$    &   $1.70{\times}10^{-4}$  & $0.65$    & $18.42$ \\
\hline
      $2705$    &   $4.71{\times}10^{-4}$  & $0.21$    & $11.04$ \\
\hline
       $3471$   &  $0.014$ & $0.07$  & $5.03$  \\

\hline
\end{tabular}

\caption{Relative errors between the fine scale solution and {\bf offline spaces} obtained from 
              using oversampled domains; Eigenvalue problem from
              Eq.~\eqref{offeig3}, $20{\times}20$ coarse mesh, harmonic snapshots, $\omega_i^+=3\times\omega_i$.}

\label{off20_parindep_At_Sp}
\end{table}
%

\begin{table}[htb]
\centering

\begin{tabular}{|c|c|c|c|c|c|}
\hline 
\multirow{2}{*}{$\text{dim}(V_{\text{off}}^+)$} &\multirow{2}{*}{$\Lambda_*$} &
\multicolumn{2}{c|}{  $\|u^{+,3471}-u^{+, \text{off}} \|$ (\%) }  \\
\cline{3-4} {}&{}&
$\hspace*{0.8cm}   L^{2}_\kappa(D)   \hspace*{0.8cm}$ &
$\hspace*{0.8cm}   H^{1}_\kappa(D)  \hspace*{0.8cm}$
\\
\hline\hline
       $1524$     &     $3.25{\times}10^{-5}$   &  $2.72$    & $35.98$  \\
\hline
      $2168$    &   $1.70{\times}10^{-4}$  & $0.69$    & $17.70$ \\
\hline
      $2705$    &   $4.71{\times}10^{-4}$  & $0.22$    & $9.82$ \\
\hline
       $3182$   &  $0.0059$ & $0.02$  & $3.21$  \\

\hline
\end{tabular}

\caption{Relative errors between the maximal dimension offline solution ($u^{+,3471}$) and 
              {\bf offline spaces} obtained from 
              using oversampled domains; Eigenvalue problem from
              Eq.~\eqref{offeig3}, $20{\times}20$ coarse mesh, harmonic snapshots, $\omega_i^{+}=3\times \omega_i$.}

\label{off20_parindep_At_Sp_vs3471}
\end{table}

\begin{table}[htb]
\centering

\begin{tabular}{|c|c|c|c|c|c|}
\hline 
\multirow{2}{*}{$\text{dim}(V_{\text{off}})$} &\multirow{2}{*}{$\Lambda_*$} &
\multicolumn{2}{c|}{  $\|u-u^{\text{off}} \|$ (\%) }  \\
\cline{3-4} {}&{}&
$\hspace*{0.8cm}   L^{2}_\kappa(D)   \hspace*{0.8cm}$ &
$\hspace*{0.8cm}   H^{1}_\kappa(D)  \hspace*{0.8cm}$
\\
\hline\hline
       $1524$ & $0.03$ &  $0.27$  & $19.47$  \\
\hline
      $2607$ & $0.06$  & $0.15$ & $12.39$ \\
\hline
      $3690$ & $0.16$  & $0.07$ & $9.40$ \\
\hline
      $7300$ & $684.14$& $0.05$ &$3.70$  \\

\hline
\end{tabular}

\caption{Relative errors between the fine scale solution and {\bf offline spaces} obtained from 
              using oversampled domains; Eigenvalue problem from
              Eq.~\eqref{offeig3}, $20{\times}20$ coarse mesh, harmonic snapshots, $\omega_i^{+}=\omega_i+1$.}

\label{off20_parindep_At_Stt1+}
\end{table}

\begin{table}[htb]
\centering

\begin{tabular}{|c|c|c|c|c|c|}
\hline 
\multirow{2}{*}{$\text{dim}(V_{\text{off}}^+)$} &\multirow{2}{*}{$\Lambda_*$} &
\multicolumn{2}{c|}{  $\|u-u^{+, \text{off}} \|$ (\%) }  \\
\cline{3-4} {}&{}&
$\hspace*{0.8cm}   L^{2}_\kappa(D)   \hspace*{0.8cm}$ &
$\hspace*{0.8cm}   H^{1}_\kappa(D)  \hspace*{0.8cm}$
\\
\hline\hline
       $1524$     &     $0.002$   &  $1.25$    & $28.54$  \\
\hline
      $2102$    &   $0.009$  & $0.22$    & $14.34$ \\
\hline
      $2607$    &   $0.014$  & $0.12$    & $8.25$ \\
\hline
       $3596$   &  $1.03{\times}10^{3}$ & $0.01$  & $2.06$  \\

\hline
\end{tabular}

\caption{Relative errors between the fine scale solution and {\bf offline spaces} obtained from 
              using oversampled domains; Eigenvalue problem from
              Eq.~\eqref{offeig3}, $20{\times}20$ coarse mesh, eigenvalue snapshots, $\omega_i^+=3\times\omega_i$.}

\label{off20_parindep_At_Sp_eigen_t5}
\end{table}

\begin{table}[htb]
\centering

\begin{tabular}{|c|c|c|c|c|c|}
\hline 
\multirow{2}{*}{$\text{dim}(V_{\text{off}}^+)$} &\multirow{2}{*}{$\Lambda_*$} &
\multicolumn{2}{c|}{  $\|u-u^{+, \text{off}} \|$ (\%) }  \\
\cline{3-4} {}&{}&
$\hspace*{0.8cm}   L^{2}_\kappa(D)   \hspace*{0.8cm}$ &
$\hspace*{0.8cm}   H^{1}_\kappa(D)  \hspace*{0.8cm}$
\\
\hline\hline
       $1524$     &     $0.23$   &  $0.75$    & $21.01$  \\
\hline
      $2168$    &   $0.82$  & $0.18$    & $13.19$ \\
\hline
      $2705$    &   $2.45$  & $0.08$    & $8.67$ \\
\hline
       $3471$   &  $117.18$ & $0.07$  & $4.44$  \\

\hline
\end{tabular}

\caption{Relative errors between the fine scale solution and {\bf offline spaces} obtained from 
              using oversampled domains with $\kappa=\frac{\kappa} {10^{4}}$ in $\omega_{i}^{+}\backslash\omega_{i}$; Eigenvalue problem from
              Eq.~\eqref{offeig3}, $20{\times}20$ coarse mesh, harmonic snapshots, $\omega_i^+=3\times\omega_i$.}

\label{off20_parindep_At_SpChoppedCoeff}
\end{table}

\begin{figure}[htb]
 \centering
 \subfigure[$N=10$, corrcoef=0.98]{\label{fig:ErvsLambdaN10}
    \includegraphics[width = 0.45\textwidth, keepaspectratio = true]{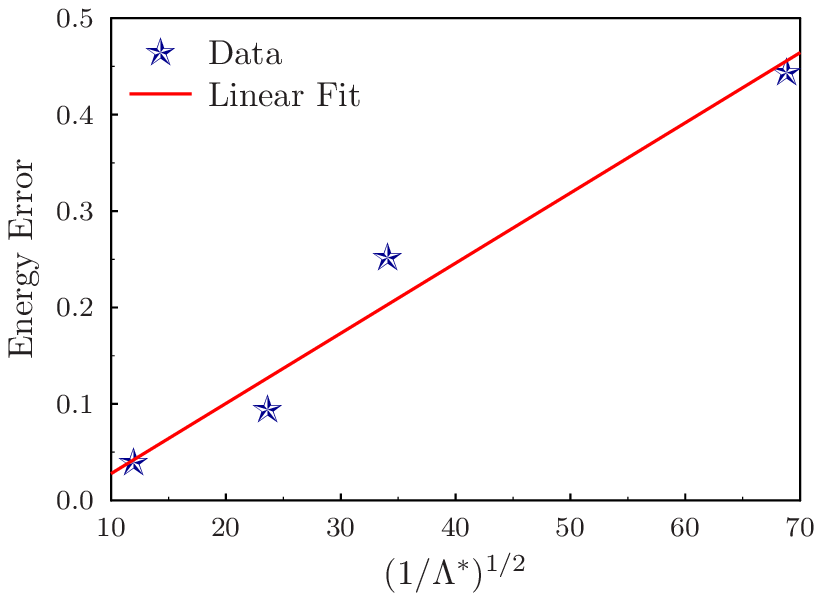}
   }
  \subfigure[$N=20$, corrcoef=0.99]{\label{fig:ErvsLambdaN20}
     \includegraphics[width = 0.45\textwidth, keepaspectratio = true]{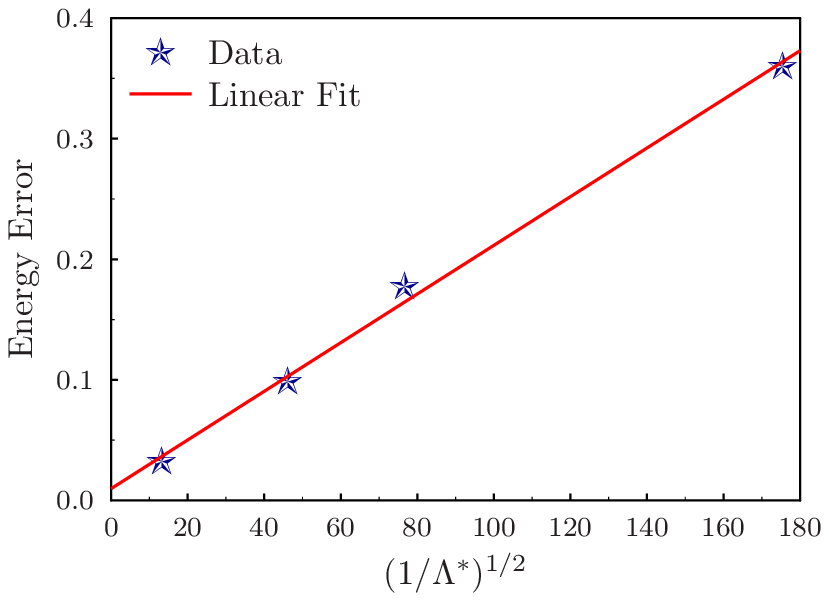}
  }
\subfigure[$N=10$, corrcoef=0.93, eigenvalue snaps]{\label{fig:ErvsLambdaN10_eigsnap}
     \includegraphics[width = 0.45\textwidth, keepaspectratio = true]{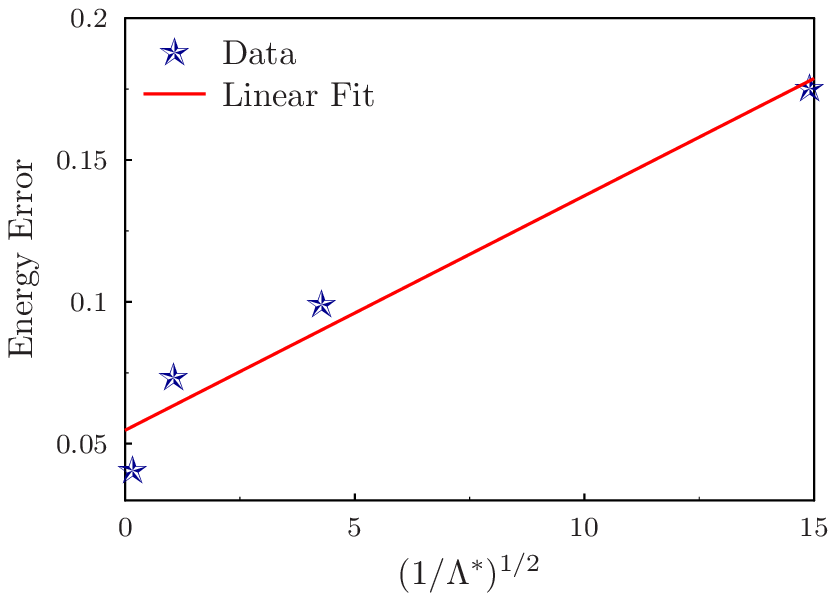}
  }
  \subfigure[$N=20$, corrcoef=0.97, eigenvalue snaps]{\label{fig:ErvsLambdaN20_eigsnap}
       \includegraphics[width = 0.45\textwidth, keepaspectratio = true]
       {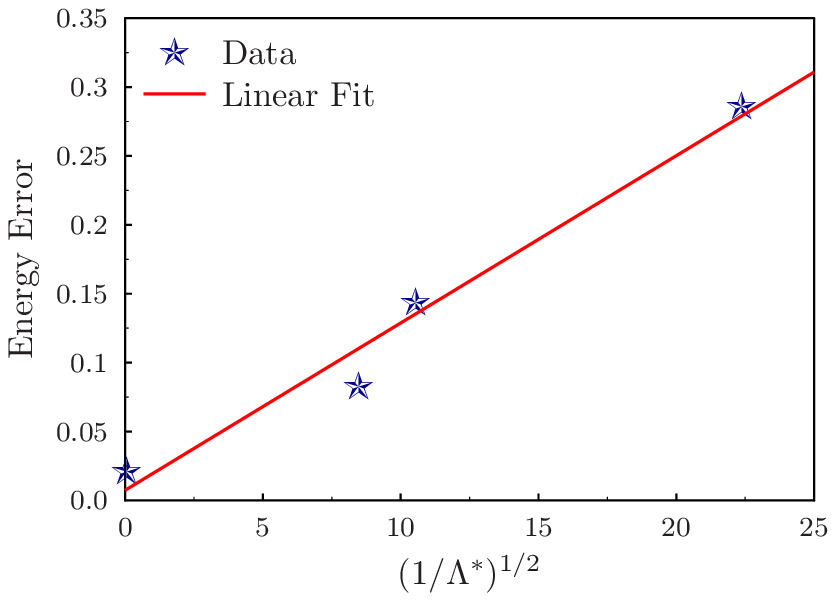}
    }
 \caption{Relation between relative energy error and $\Lambda_*$ for Tables \ref{off10_parindep_At_Sp_vs909}, \ref{off20_parindep_At_Sp_vs3471}, \ref{off10eigsnap_parindep_At_Sp} and \ref{off20_parindep_At_Sp_eigen_t5} respectively.}
\label{fig:ErvsLambda}
\end{figure}

Finally, we plot the energy error against $(1/\Lambda_*)^{1\over 2}$ for 
$10 \times 10$ and $20 \times 20$ cases in Figs.~\ref{fig:ErvsLambda}. 
The correlation between
the errors and  $1/\Lambda_*$ is over $0.93$ when we consider $10 \times 10$
 mesh case 
(as in Figs.~\ref{fig:ErvsLambdaN10} and \ref{fig:ErvsLambdaN10_eigsnap}). 
In Figs.~~\ref{fig:ErvsLambdaN20} and
\ref{fig:ErvsLambdaN20_eigsnap}, we depict the
relative errors corresponding to Tables 
\ref{off20_parindep_At_Sp_vs3471} and
\ref{off20_parindep_At_Sp_eigen_t5}.
In this case, we also observe a good agreement and
 the correlation to be over $0.98$.


\subsubsection{Parameter-independent case using multiple spectral problems}
\label{sec:num_res_multiple}

In this section, we study the use of multiple spectral problems
as described in Section \ref{sec:multiple}. In particular,
we use only two spectral problems in $\omega$ and $\omega^{+}$.
The results are presented in
Tables \ref{off10_parindep_At_St_withTol_subss} and 
\ref{off10_parindep_At_St_withTol_subss_off-snap}. 
As the convergence theory indicates, for the same 
eigenvalue threshold, one can expect the quadratic decay in the 
convergence rate with a constant that is described 
in Section \ref{appendix:convergence_multiple}.
In Table \ref{off10_parindep_At_St_withTol_subss},
we compare the offline solution and the fine grid solution,
while in Table \ref{off10_parindep_At_St_withTol_subss_off-snap},
we compare the offline solution and the snapshot solution.
In both cases, we observe that the square of the error resulting from
 a single spectral problem correlates well to the case
corresponding to multiple spectral problems. This behavior deteriorates
when the space dimension is large due to irreducible error.
For this set of numerical results, we observe that the coarse space
dimension resulting from multiple spectral problems is large compared
to the case when a single spectral problem is used. However, we note that
our convergence result does not contain any information about
the dimension of the coarse space, but only about an eigenvalue
threshold. On the other hand, our convergence analysis suggests
that the coarse space needs to include an approximation in both
$\omega$ and $\omega^{+}$. The eigenvectors of Eq.~\eqref{offeig4}
may be represented using  the eigenvectors of Eq.~\eqref{offeig1},
and thus one can use the respective eigenvectors to complement each other.
Our numerical results show that by combining eigenvectors
of  Eq.~\eqref{offeig1} and  Eq.~\eqref{offeig4}, one can
achieve better convergence compared to only using  Eq.~\eqref{offeig1}
in our pre-asymptotic numerical simulations.

%
%
\begin{table}[htb]
\centering
{\begin{tabular}{|c|c|c|c|}

\hline 
{$\text{dim}(V_{\text{off}})$} &{$\Lambda^*$} &
 \scriptsize{${H^{1}_\kappa(D)}$ (\%) (Eqs.~\eqref{offeig4}, ~\eqref{offeig1}) }&
 \scriptsize{ ${H^{1}_\kappa(D)}$ (\%) (Eq.~\eqref{offeig1}) }  \\
\hline\hline
       $791(618)$  &   $60.52(147.82)$ ($\text{tol}= 60(100)$)     & $31.47$&   $38.34$  \\
\hline
      $1172(733)$  &   $401.56(1.01{\times}10^3)$ ($\text{tol}=400(1000)$)   & $14.72$    & $24.42$ \\
\hline
      $2054(1568)$    &   $1.00(5.05){\times}10^3$ ($\text{tol}=1000(5000)$)    & $8.28$    & $10.73$ \\
\hline
\end{tabular}
}
\caption{Relative errors between the fine scale solution and {\bf offline spaces} for local spectral problems using a single (Eq.~ \eqref{offeig1}) 
and multiple eigenvalue problems (Eq.~\eqref{offeig4} and \eqref{offeig1}).
 $20{\times}20$ coarse mesh, harmonic snapshots, $\omega_i^+=3\times\omega_i$.}
\label{off10_parindep_At_St_withTol_subss}
\end{table}

\normalsize

%
%
\begin{table}[htb]
\centering

{\begin{tabular}{|c|c|c|c|}

\hline 
{$\text{dim}(V_{\text{off}})$} &{$\Lambda^*$} &
 \scriptsize{${H^{1}_\kappa(D)}$ (\%) (Eqs.~\eqref{offeig4}, ~\eqref{offeig1}) }&
 \scriptsize{ ${H^{1}_\kappa(D)}$ (\%) (Eq.~\eqref{offeig1}) }  \\
\hline\hline
       
$791(618)$  &   $60.52(147.82)$ ($\text{tol}= 60(100)$) & $31.26$ & $38.16$  \\
\hline
      $1172(733)$  &   $401.56(1.01{\times}10^3)$ ($\text{tol}=400(1000)$) & $14.31$ & $24.16$ \\
\hline
      $2054(1568)$    &   $1.00(5.05){\times}10^3$ ($\text{tol}=1000(5000)$) & $7.53$ & $10.16$ \\
\hline
\end{tabular}
}
\caption{Relative errors between the snapshot solution and {\bf offline spaces} for local spectral problems using a single (Eq.~ \eqref{offeig1}) 
and multiple eigenvalue problems (Eq.~\eqref{offeig4} and \eqref{offeig1}).
 $20{\times}20$ coarse mesh, harmonic snapshots, $\omega_i^+=3\times\omega_i$.}
\label{off10_parindep_At_St_withTol_subss_off-snap}
\end{table}

\normalsize
%

\subsection{Parameter-dependent case}\label{sec:pdcase}

For the next set of numerical results, we consider a parameter-dependent
example where
\[
\kappa(x;\mu)=\mu_1 \kappa_1(x) + \mu_2\kappa_2(x),
\]
where $\kappa_1(x)$ and $\kappa_2(x)$ are depicted in Fig.~\ref{perm}.
For the numerical examples, we consider a snapshot space that consists
of solving local eigenvalue problem described by Eq.~\eqref{offeig1} in
$\omega_i^{+}=3\times \omega_i$ for 9 selected values of $\mu=(\mu_1,\mu_2)$.
By choosing 20 dominant eigenvectors for each of 9 selected values
of $\mu$ and ensuring linear independence, we form the space of snapshots. In this space of snapshots,
we use the operator averaged over $\mu$ to construct the offline space.
In particular, we consider $3$ choices for offline eigenvalue problems
that are given by
Eq.~\eqref{offeig1}, Eq.~\eqref{offeig2}, and Eq.~\eqref{offeig3}.
These local eigenvalue problems are used  to construct the offline
space. Furthermore, we use the same eigenvalue problems for an online
value of the parameter, to construct an online space which
is a subspace of the offline space by solving the local eigenvalue
problem Eq.~\eqref{oneig1}, Eq.~\eqref{oneig2}, and Eq.~\eqref{oneig3}.
The results are presented in Tables \ref{on_pardep_At_Mt},
\ref{on_pardep_Ap_At}, \ref{on_pardep_At_Sp}, respectively.
We see from these results
that the GMsFEM converges in all the cases considered above. The best convergence  among the three 
choices it found in Table \ref{on_pardep_At_Sp}.

\begin{table}[htb]
\centering
\begin{tabular}{|c|c|c|c|c|c|}
\hline \multirow{2}{*}{$\text{dim}(V_{\text{on}})$}
&\multirow{2}{*}{$\Lambda_*$}&
\multicolumn{2}{c|}{  $\|u-u^{\text{on}} \|$ (\%) }&\multicolumn{2}{c|}{ $\|u^{\text{on}}-u^{\text{off}} \|$ (\%)}\\
\cline{3-6} {}&{}&
$\hspace*{0.8cm}   L^{2}_\kappa(D)   \hspace*{0.8cm}$ &
$\hspace*{0.8cm}   H^{1}_\kappa(D)  \hspace*{0.8cm}$&
$\hspace*{0.8cm}   L^{2}_\kappa(D)   \hspace*{0.8cm}$&
$\hspace*{0.8cm}   H^{1}_\kappa(D)  \hspace*{0.8cm}$
\\
\hline\hline
       $728$     &    $537.1$   & $0.26$    &$11.20$   &  $0.18$   &   $9.50$ \\
\hline
      $907$     &   $1.05{\times}10^3$    & $0.18$    &$8.62$  & $0.09$   &   $6.27$ \\
\hline
       $1037$    &  $1.93{\times}10^3$    & $0.17$  &$8.22$    & $0.09$   &   $5.72$\\
\hline
      $1248$    &  $3.22{\times}10^3$        &  $0.11$    &$6.42$ &     $0.02$   &$2.55$ \\
\hline
      $1378$    &  ---        &  $0.10$    &$5.89$ &     $0.00$   &$0.00$ \\
\hline
\end{tabular}

\caption{Relative errors between the fine scale solution (and offline) and {\bf online spaces} obtained from using non-oversampled domains; Eigenvalue problem from Eq.~\eqref{oneig1}, $10{\times}10$ coarse mesh, eigenvalue snapshots, $\omega_i^{+}=3\times \omega_i$.}
\label{on_pardep_At_Mt}
\end{table}

\begin{table}[htb]
\centering
\begin{tabular}{|c|c|c|c|c|c|}
\hline \multirow{2}{*}{$\text{dim}(V_{\text{on}}^+)$}
&\multirow{2}{*}{$\Lambda_*$}&
\multicolumn{2}{c|}{  $\|u-u^{+,\text{on}} \|$ (\%) }&\multicolumn{2}{c|}{ $\|u^{+,\text{on}}-u^{+,\text{off}} \|$ (\%)}\\
\cline{3-6} {}&{}&
$\hspace*{0.8cm}   L^{2}_\kappa(D)   \hspace*{0.8cm}$ &
$\hspace*{0.8cm}   H^{1}_\kappa(D)  \hspace*{0.8cm}$&
$\hspace*{0.8cm}   L^{2}_\kappa(D)   \hspace*{0.8cm}$&
$\hspace*{0.8cm}   H^{1}_\kappa(D)  \hspace*{0.8cm}$
\\
\hline\hline
       $826$     &     \hspace*{0.2cm} $11.47$   \hspace*{0.2cm}  & $0.27$    &$10.93$   &  $0.18$   &   $9.19$ \\
\hline
      $988$     &   $40.84$    & $0.17$    &$8.37$  & $0.09$   &   $5.93$ \\
\hline
       $1133$    &  $65.10$    & $0.10$  &$6.23$    & $0.01$   &   $2.02$\\
\hline
      $1280$    &  $350.0$        &  $0.10$    &$6.03$ &     $0.007$   &$1.27$ \\
\hline
      $1378$    &  ---        &  $0.10$    &$5.89$ &     $0.00$   &$0.00$ \\

\hline
\end{tabular}

\caption{Relative errors between the fine scale (and offline) solution and {\bf online spaces} obtained from using oversampled domains; Eigenvalue problem from Eq.~\eqref{oneig2}, $10{\times}10$ coarse mesh, eigenvalue snapshots, $\omega_i^{+}=3\times \omega_i$. }
\label{on_pardep_Ap_At}
\end{table}

\begin{table}[htb]
\centering
\begin{tabular}{|c|c|c|c|c|c|}
\hline \multirow{2}{*}{$\text{dim}(V_{\text{on}}^+)$}
&\multirow{2}{*}{$\Lambda_*$}&
\multicolumn{2}{c|}{  $\|u-u^{+,\text{on}} \|$ (\%) }&\multicolumn{2}{c|}{ $\|u^{+,\text{on}}-u^{+,\text{off}} \|$ (\%)}\\
\cline{3-6} {}&{}&
$\hspace*{0.8cm}   L^{2}_\kappa(D)   \hspace*{0.8cm}$ &
$\hspace*{0.8cm}   H^{1}_\kappa(D)  \hspace*{0.8cm}$&
$\hspace*{0.8cm}   L^{2}_\kappa(D)   \hspace*{0.8cm}$&
$\hspace*{0.8cm}   H^{1}_\kappa(D)  \hspace*{0.8cm}$
\\
\hline\hline
       $790$     &     $4.33{\times}10^{-4}$ & $0.22$    &$10.37$   &  $0.25$   &   $8.51$ \\
\hline
      $888$     &   $0.0032$    & $0.09$    &$7.90$  & $0.06$   &   $5.25$ \\
\hline
       $1133$    &  $3.27$    & $0.10$  &$6.29$    & $0.03$   &   $2.18$\\
\hline
      $1280$    &  $154.3$        &  $0.10$    &$5.96$ &     $0.004$   &$0.87$ \\
\hline
      $1378$    &  ---        &  $0.10$    &$5.89$ &     $0.00$   &$0.00$ \\

\hline
\end{tabular}

\caption{Relative errors between the fine scale (and offline) solution and {\bf online spaces} obtained from using oversampled domains; Eigenvalue problem from Eq.~\eqref{oneig3}, $10{\times}10$ coarse mesh, eigenvalue snapshots, $\omega_i^{+}=3\times \omega_i$. }
\label{on_pardep_At_Sp}
\end{table}

\section{Convergence studies and discussions.
 Parameter-independent case}
\label{appendix:convergence}

\subsection{Offline space using a single spectral problem}
\label{appendix:convergence_single}

We define $I^{\omega_i}u$ and
$I^{\omega_i^{+}}u$ as an interpolation of $u$ in
$\omega_i$ and $\omega_i^{+}$ that will be chosen.
Because the snapshot functions are defined in $\omega_i^{+}$,
$I^{\omega_i}u=I^{\omega_i^{+}}u$ in $\omega_i$.
 We have
\begin{equation}
\label{eq:residual}
\begin{split}
-\mbox{div}(\kappa(x) \nabla (u-I^{\omega_i}u))=g\ \text{in}\ \omega_i,\\
-\mbox{div}(\kappa(x) \nabla (u-I^{\omega_i^{+}}u))=g\ \text{in}\ \omega_i^{+},
\end{split}
\end{equation}
where $g$ is the source term, $g=f+\mbox{div}(\kappa(x) \nabla I^{\omega_i}u)$.
Furthermore, we define $\chi_i$ and $\chi_i^{+}$ as partition
of unity functions subordinated to $\omega_i$ and
$\omega_i^{+}$. In particular, we can assume that 
 $\chi_i$ and $\chi_i^{+}$  are standard bilinear partition
of unity functions for a rectangular partition.
In general, we do not require $\chi_i^{+}$ to be a partition
of unity function; however, we require $\chi_i^{+}$ to be zero
on $\partial \omega_i^{+}$ and 
\[
|\nabla \chi_i|^2\preceq |\nabla \chi_i^{+}|^2.
\]
This is satisfied for bilinear functions.

Multiplying both sides
of (\ref{eq:residual}) by $\chi_i^2(u-I^{\omega_i}u)$
(or $(\chi_i^{+})^2(u-I^{\omega_i^{+}}u)$ for the equation in $\omega_i^{+}$), 
integrating by parts
 and re-arranging the terms,
we have
\begin{eqnarray*}
&&\int_{\omega_i} \kappa \chi_i^2 |\nabla(u-I^{\omega_i} u)|^2\\
& \leq&
 {1\over C}\int_{\omega_i} \kappa  |\nabla \chi_i|^2|(u-I^{\omega_i}u)|^2 +
 C \int_{\omega_i}  \kappa \chi_i^{2} |\nabla(u-I^{\omega_i} u)|^{2}+
 |\int_{\omega_i} g \chi_i^2 (u-I^{\omega_i}u)|,
\end{eqnarray*}
where $C<1$ is independent of contrast.
From here, we get  
\begin{equation}\label{eq:CinOmega}
\begin{split}
\int_{\omega_i} \kappa \chi_i^2 |\nabla(u-I^{\omega_i} u)|^2 \preceq
 \int_{\omega_i} \kappa  |\nabla \chi_i|^2|(u-I^{\omega_i}u)|^2 +
 |\int_{\omega_i} g \chi_i^2 (u-I^{\omega_i}u)|.
\end{split}
\end{equation}
Similarly,
\begin{equation}\label{eq:CinOmega+}
\begin{split}
\int_{\omega_i^{+}} \kappa |\chi_i^{+}|^2 |\nabla(u-I^{\omega_i^{+}} u)|^2 \preceq
 \int_{\omega_i^{+}} \kappa  |\nabla \chi_i^{+}|^2|(u-I^{\omega_i^{+}}u)|^2 +
 |\int_{\omega_i^{+}} g |\chi_i^{+}|^2 (u-I^{\omega_i^{+}}u)|.
\end{split}
\end{equation}
Next, 
taking into account that MsFEM solution, $u_H$,
 provides a minimal energy error,
we have
\begin{eqnarray}
\int_D \kappa  |\nabla(u-u_H)|^2
&\preceq& \int_D \kappa  |\nabla(\sum_i \chi_i(u-I^{\omega_i}u))|^2\nonumber\\ 
&\preceq&
\sum_i \int_{\omega_i} \kappa |\nabla \chi_i|^2 |u-I^{\omega_i}u|^2 +
\sum_i \int_{\omega_i} \kappa | \chi_i|^2 |\nabla(u-I^{\omega_i}u)|^2.
\end{eqnarray}
Combining this inequality with~\eqref{eq:CinOmega}, we obtain,
\begin{eqnarray}
\int_D \kappa  |\nabla(u-u_H)|^2
&\preceq&
\sum_i \int_{\omega_i} \kappa |\nabla \chi_i|^2 |u-I^{\omega_i}u|^2 +\left(
\sum_i \int_{\omega_i} \kappa | \chi_i|^2 |\nabla(u-I^{\omega_i}u)|^2\right)\nonumber\\
&\preceq& \sum_i \int_{\omega_i} \kappa |\nabla \chi_i|^2 |u-I^{\omega_i}u|^2
+
\left( { \int_{\omega_i} \kappa  |\nabla \chi_i|^2|(u-I^{\omega_i}u)|^2}\right. \nonumber\\
&&+\left.
\sum_i|\int_{\omega_i} g |\chi_i|^2(u-I^{\omega_i}u)|\right) \label{eq:16}\\
&\preceq& \sum_i \int_{\omega_i} \kappa |\nabla \chi_i|^2 |u-I^{\omega_i}u|^2
+
\sum_i|\int_{\omega_i} g |\chi_i|^2(u-I^{\omega_i}u)| \nonumber
\end{eqnarray}
{Next, we concentrate in deriving a bound for the first term on the right hand side of the last inequality above}.

Note that $I^{\omega_i^{+}}u=I^{\omega_i}u$ in $\omega_i$. 
Next, we define the interpolant $I^{\omega_i^{+}}u$ using 
the modes for the eigenvalue problem Eq.~\eqref{offeig2}
that correspond to the eigenvalues
$\lambda^{\omega_i}_{1}, \cdots,\lambda^{\omega_i}_{L_i}$. Then, we have
\begin{equation}\label{eq:fact}
\int_{\omega_i^{+}}\kappa |\nabla \chi_i^{+}|^2 (u-I^{\omega_i^{+}}u)^2 \preceq  {1\over \lambda^{\omega_i}_{L_i+1}}
\int_{\omega_i} \kappa  |\nabla (u-I^{\omega_i}u)|^2,
\end{equation} 
{which is easily deduced 
from the corresponding eigenvalue problem and the definition of the interpolation $I^{\omega_i}u$},
we have
\begin{equation}\label{eq:bound1}
\sum_i \int_{\omega_i} \kappa |\nabla \chi_i|^2 |u-I^{\omega_i}u|^2\preceq
\sum_i \int_{\omega_i^{+}} \kappa |\nabla \chi_i^{+}|^2 |u-I^{\omega_i^{+}}u|^2 
\preceq \sum_i {1\over \lambda^{\omega_i}_{L_i+1}}
\int_{\omega_i} \kappa  |\nabla (u-I^{\omega_i}u)|^2.
\end{equation}
Note also that we can bound the last term above by
\begin{eqnarray*}
&&\sum_i {1\over \lambda^{\omega_i}_{L_i+1}}
\int_{\omega_i} \kappa  |\nabla (u-I^{\omega_i}u)|^2\preceq
{\sum_i
 {1\over \lambda^{\omega_i}_{L_i+1}}
\int_{\omega_i^{+}} \kappa  |\chi_i^{+}|^2 |\nabla (u-I^{\omega_i^{+}}u)|^2}\\
&\preceq& \sum_i {1 \over \lambda^{\omega_i}_{L_i+1}} \int_{\omega_i^{+}}
 \kappa |\nabla \chi_i^{+}|^2 |u-I^{\omega_i^{+}}u|^2 + \sum_i {1\over
\lambda^{\omega_i}_{L_i+1}} |\int_{\omega_i^{+}} g |\chi_i^{+}|^2 (u-I^{\omega_i^{+}} u)|  \quad
\mbox{{ (by~\eqref{eq:CinOmega+})}}\\
&\preceq&
{1\over\Lambda_*}
\left( \sum_i  \int_{\omega_i^{+}} \kappa |\nabla \chi_i^{+}|^2 |u-I^{\omega_i^{+}}u|^2 + \sum_i  |\int_{\omega_i^{+}} g |\chi_i^{+}|^2 (u-I^{\omega_i^{+}} u)| \right)\\
&\preceq&
{1\over\Lambda_*}
\left(   \sum_i {1\over \lambda^{\omega_i}_{L_i+1}}
\int_{\omega_i} \kappa  |\nabla (u-I^{\omega_i}u)|^2 + \sum_i  |\int_{\omega_i^{+}} g |\chi_i^{+}|^2 (u-I^{\omega_i^{+}} u)| \right) 
, \quad \mbox{{ (by~\eqref{eq:fact})}}\nonumber
\end{eqnarray*}
where {we have defined} $\Lambda_*=\min_{\omega_i}
\lambda_{L_i+1}^{\omega_i}$.
Thus, {summarizing the last set of inequalities we obtain},
\[
\sum_i {1\over \lambda^{\omega_i}_{L_i+1}}
\int_{\omega_i} \kappa  |\nabla (u-I^{\omega_i}u)|^2\preceq
{1\over\Lambda_*}
\left(   \sum_i {1\over \lambda^{\omega_i}_{L_i+1}}
\int_{\omega_i} \kappa  |\nabla (u-I^{\omega_i}u)|^2 + \sum_i  |\int_{\omega_i^{+}} g |\chi_i^{+}|^2 (u-I^{\omega_i^{+}} u)| \right).
\]

Applying this inequality $n$ times in the estimate for
$\sum_i \int_{\omega_i} \kappa |\nabla \chi_i|^2 |u-I^{\omega_i}u|^2$ 
{in Eq.~\eqref{eq:bound1}}, we get
\begin{eqnarray*}
&&\sum_i \int_{\omega_i^{+}} \kappa |\nabla \chi_i^{+}|^2 |u-I^{\omega_i^{+}}u|^2\preceq
\sum_i {1\over \lambda^{\omega_i}_{L_i+1}}
\int_{\omega_i} \kappa  |\nabla (u-I^{\omega_i}u)|^2\\
 &\preceq&\left({1\over \Lambda_*}\right)^n \sum_i \int_{\omega_i}  {1\over \lambda^{\omega_i}_{L_i+1}}\kappa
 |\nabla(u-I^{\omega_i}u)|^2 + \sum_{j=1}^{n}
\left({1\over \Lambda_*}\right)^j \sum_i |
\int_{\omega_i^{+}} g |\chi_i^{+}|^2 (u-I^{\omega_i^{+}} u)|\\
 &\preceq&\left({1\over \Lambda_*}\right)^{n+1} \sum_i \int_{\omega_i} \kappa
 |\nabla(u-I^{\omega_i}u)|^2 + \left({\Lambda_*}\right)^n
\left({1-\Lambda_*^{-n}\over \Lambda_*-1}\right)
\sum_i \int_{\omega_i^{+}} (| \kappa|| \nabla \chi_i^{+}|^2)^{-1} g^2.
\end{eqnarray*}
Considering
$\sum_i \int_{\omega_i} \kappa |\nabla(u-I^{\omega_i}u)|^2 \preceq
\int_{D} \kappa  |\nabla u|^2$,
we have the following convergence rate for GMsFEM,
\begin{equation}
\label{eq:rate1}
\begin{split}
\int_D \kappa  |\nabla(u-u_H)|^2 \preceq
\frac{1}{ \Lambda_*^{n+1} }\int_{D} \kappa |\nabla u|^2 +
\left( \left({\Lambda_*}\right)^{n} \left({1-\Lambda_*^{-n}\over \Lambda_*-1}\right) +1 \right) R,
\end{split}
\end{equation}
where $R=\sum_i \int_{\omega_i^{+}}(| \kappa||\nabla \chi_i^{+}|^2)^{-1} g^2$. 
For right hand sides {with} $g\preceq 1$
it can be shown that
$\int_{\omega_i^{+}} \left ( \kappa |\nabla \chi_i^{+}|^2\right)^{-1} g^2\preceq
H^2$.

With this assumption, { we have the  convergence result},
\begin{equation}
\label{eq:rate2}
\begin{split}
\int_D \kappa  |\nabla(u-u_H)|^2 \preceq
\frac{1}{ \Lambda_*^{n+1}} \int_{D} \kappa |\nabla u|^2 +
\left( \left({\Lambda_*}\right)^n \left({1-\Lambda_*^{-n}\over \Lambda_*-1}\right) +1 \right)H^2\int_{D} |1|^2.
\end{split}
\end{equation}
Choosing $ \Lambda_*$ sufficiently large (larger than $1$)
and $n= 1-{\log(H)\over \log{\Lambda_*}}$ (in each $\omega_i$),
we obtain
\begin{equation}
\label{eq:rate13}
\begin{split}
\int_D \kappa  |\nabla(u-u_H)|^2 \preceq
  \left( {H\over\Lambda_*}  \right) \int_D \kappa |\nabla u|^2 + {H\over\Lambda_*}.
\end{split}
\end{equation}
Collecting the results above, we have 
\begin{theorem}
If $\Lambda_*\geq 1$ and $\int_D \kappa^{-1}g^2\preceq 1$, then
\[
\int_D \kappa  |\nabla(u-u_H)|^2 \preceq
 \left( {H\over\Lambda_*}  \right) \int_D \kappa |\nabla u|^2 + {H\over\Lambda_*}.
\]
\end{theorem}

Next, we comment on the estimate on $g$. We
consider the snapshot space generated
by Eq.~\eqref{snaplinalg}.
Because $|\nabla \chi_i^{+}|^2\preceq H^{-2}$, we have
$R\preceq H^2\int_{\omega_i^{+}}\kappa^{-1} g^2$. We assume that
$I^{\omega_i^{+}}u=\sum_{l=1}^L c_l\Psi_l$, where $\Psi_l$ are eigenvectors
Eq.~\eqref{offeig3}. Each eigenvector  $\Psi_l$ is spanned by
eigenvectors of Eq.~\eqref{snaplinalg}, i.e., $\Psi_l=d_{lm}\psi_m^{+}$.
Then,
$g=f+\mbox{div}(\kappa(x) \nabla I^{\omega_i}u)=
f-\sum_{m=1}^{M_{\text{snap}}} \sum_{l=1}^L c_ld_{lm}\lambda_m\kappa \psi_m^{+}=f-\sum_m d_m^*\lambda_m\kappa \psi_m^{+}$, where $d_{m,L}^*= \sum_{l=1}^L c_ld_{lm}$,
$\lambda_m$ are eigenvalues in Eq.~\eqref{snaplinalg}, $M_{\text{snap}}$
is the number of snapshots, and $L$ is the number of modes selected in the
offline stage.
Due to orthogonality of
$\psi_m^{+}$, it can be shown that 
$\int_{\omega_i^{+}}\kappa^{-1} g^2\preceq 1 + \sum_m (d_{m,L}^*)^2 \lambda_m^2$,
provided $1\preceq \kappa$.
On the other hand, $\int_{\omega_i^+}\kappa|\nabla I^{\omega_i^+}u|^2 =
\sum_m (d_{m,L}^*)^2 \lambda_m$. Thus, 
\[
\int_{\omega_i^{+}}\kappa^{-1} g^2\preceq 1 +\Lambda_*^{\text{snap}} \int_{\omega_i^+}\kappa|\nabla I^{\omega_i^+}u|^2. 
\]

\subsection{Offline space using multiple spectral problems and the relation to oversampled spectral problems}
\label{appendix:convergence_multiple}

In this subsection we briefly consider the convergence analysis 
for the offline space
proposed in Section \ref{sec:multiple}. 
For this analysis we use the harmonic snapshot space to avoid any residual error, though
the derivation can be extended to other scenarios. This derivation uses the proof of
Theorem 3.3
from the work of  Babu{\v{s}}ka and Lipton \cite{bl11} and we extend it to
a high-contrast case.
We start with the inequality in 
Eq.~\eqref{eq:16}. 
We define two interpolants $I^{\omega_i^{+}}u$ and
 $I^{\omega_i}u$ by choosing the dominant modes through considering
\begin{equation}\label{eq:fact1}
\int_{\omega_i^{+}}\kappa |\nabla \chi_i^{+}|^2 (u-I^{\omega_i^{+}}u)^2 \preceq  {1\over \lambda^{\omega_i^{+}}_{L_i+1}}
\int_{\omega_i^{+}} \kappa  |\nabla (u-I^{\omega_i^{+}}u)|^2
\preceq  {1\over \lambda^{\omega_i^{+}}_{L_i+1}}
\int_{\omega_i^{+}} \kappa  |\nabla u|^2
\end{equation} 
\begin{equation}\label{eq:fact2}
\int_{\omega_i}\kappa |\nabla \chi_i|^2 (\widehat{u}-I^{\omega_i}\widehat{u})^2 \preceq  {1\over \lambda^{\omega_i}_{L_i+1}}
\int_{\omega_i} \kappa  |\nabla (\widehat{u}-I^{\omega_i}\widehat{u})|^2\preceq  {1\over \lambda^{\omega_i}_{L_i+1}}
\int_{\omega_i} \kappa  |\nabla \widehat{u}|^2,
\end{equation} 
where $\widehat{u}= u-I^{\omega_i^{+}}u$,  ${1\over \lambda^{\omega_i^{+}}_{L_i+1}}$ and
${1\over \lambda^{\omega_i}_{L_i+1}}$ are sufficiently small. We choose
the interpolant to be $I^{\omega_i^{+}}u + I^{\omega_i}\widehat{u}$.
Thus, $u-(I^{\omega_i^{+}}u + I^{\omega_i}\widehat{u})=\widehat{u}-I^{\omega_i}\widehat{u}$.
Then, we have
\begin{eqnarray}
&&\sum_i \int_{\omega_i} \kappa |\nabla \chi_i|^2 |\widehat{u}-I^{\omega_i}\widehat{u}|^2
\preceq \sum_i {1\over \lambda^{\omega_i}_{L_i+1}}
\int_{\omega_i} \kappa  |\nabla (\widehat{u}-I^{\omega_i}\widehat{u})|^2 \nonumber\\
&\preceq& {\sum_i
 {1\over \lambda^{\omega_i}_{L_i+1}}
\int_{\omega_i^{+}} \kappa  |\chi_i^{+}|^2 |\nabla (\widehat{u}-I^{\omega_i^{+}}\widehat{u})|^2}
\preceq \sum_i {1 \over \lambda^{\omega_i}_{L_i+1}} \int_{\omega_i^{+}}
 \kappa |\nabla \chi_i^{+}|^2 |\widehat{u}-I^{\omega_i^{+}}\widehat{u}|^2   \nonumber\\
&\preceq&
{1\over\Lambda_*}
\left( \sum_i  \int_{\omega_i^{+}} \kappa |\nabla \chi_i^{+}|^2 |\widehat{u}|^2  \right)\nonumber\\
&\preceq&
{1\over\Lambda_*}
\left(   \sum_i {1\over \lambda^{\omega_i^{+}}_{L_i+1}}
\int_{\omega_i^{+}} \kappa  |\nabla (u-I^{\omega_i}u)|^2 \right) 
\preceq {1\over\Lambda_*} {1\over\Lambda_*^{+}}\int_{D} \kappa  |\nabla u|^2,
\label{eq:bound11}
\end{eqnarray}
where $\Lambda_*=\min_{\omega_i}
\lambda_{L_i+1}^{\omega_i}$ and $\Lambda_*^{+}=\min_{\omega_i}
\lambda_{L_i+1}^{\omega_i^{+}}$. Thus, choosing $\Lambda_*$
and $\Lambda_*^{+}$ to be sufficiently large, the convergence rate
can be improved. In particular, 
the final estimates involve the product of the  convergence rates with individual 
spaces.

The above results can be summarized in the
following way.
If Eq.~\eqref{eq:fact1} and Eq.~\eqref{eq:fact2} can be satisfied by
choosing appropriate interpolants in $\omega$ and
$\omega^{+}$, then 
\[
\int_D \kappa  |\nabla(u-u_H)|^2 \preceq {1\over\Lambda_*} {1\over\Lambda_*^{+}}\int_{D} \kappa  |\nabla u|^2.
\]
This result can easily be extended to use multiple eigenvalue problems (instead of two eigenvalue 
problems).

\begin{remark}
In the above proof, we rely on the estimates that bound the $L^2_\kappa$-norm
via the $H^1_\kappa$-norm in $\omega$ and $\omega^{+}$. In addition,
we use an inequality that bounds $H^1_\kappa(\omega)$ by
$L^2_\kappa(\omega^{+})$ norms based on PDE estimates.
The latter can be replaced by a third eigenvalue problem,
$A^{+,\text{off}}\Psi_k^{\text{off}} = \lambda_k^{\text{off}} S^{\text{off}} \Psi_k^{\text{off}}$, (cf. Eq.~\eqref{offeig2}),
and one can select its important modes (corresponding to largest 
eigenvalues) to 
reduce the constant relating $H^1_\kappa(\omega)$ to $L^2_\kappa(\omega^{+})$
norms.
 We have implemented
this procedure and observed a slight improvement,
at the cost of additional basis functions (similar to the numerical 
results presented in Section \ref{sec:num_res_multiple}).
\end{remark}

Next, we discuss the relation of using multiple spectral problems 
to the eigenvalue problems discussed earlier. Because of the 
optimality of local spectral spaces (cf. \cite{bl11} and see below), 
we know that the spectral problems that are similar to 
those described by Eqs.~\eqref{offeig2} and ~\eqref{offeig3}
will provide a better convergence rate compared
to those using multiple spectral problems. For the two spectral
problems described above, one can equivalently
use (cf. Eqs.~\eqref{offeig2} and ~\eqref{offeig3}) the following
local spectral problem
\begin{equation}
\label{offeig5}
A^{+,\text{off}}\Psi_k^{\text{off}}= \lambda_k^{\text{off}} S^{\text{off}} \Psi_k^{\text{off}}
\end{equation}
for construction of the offline spaces. Indeed, if 
Eq.~\eqref{offeig5} is used, we can show that (cf. Eq.~\eqref{eq:bound11})
\begin{equation}\label{eq:bound13}
\begin{split}
\sum_i \int_{\omega_i} \kappa |\nabla \chi_i|^2 |u-I^{\omega_i}u|^2
\preceq
\sum_i {1\over \widetilde{\lambda}^{\omega_i^{+}}_{L_i+1}}
\int_{\omega_i^{+}} \kappa  |\nabla (u-I^{\omega_i}u)|^2
\preceq {1\over\widetilde{\Lambda}_*} \int_{D} \kappa  |\nabla u|^2,
\end{split}
\end{equation}
where $\widetilde{\Lambda}_*=\min_{\omega_i}
\widetilde{\lambda}_{L_i+1}^{\omega_i}$. 
On the other hand, the local spectral problem
Eq.~\eqref{offeig5} provides
an optimal subspace in the following sense.
For a fixed $L_i$ dimensional subspace in $\omega_i$, a space that provides
the smallest 
\[
\max_{u}\min_{u_0} {(u-u_0)^T S^{\text{off}} (u-u_0)\over (u-u_0)^T A^{+,\text{off}} (u-u_0)}
\]
is given by the the span of the
smallest (in terms of corresponding eigenvalues)
 $L_i$ eigenvectors of Eq.~\eqref{offeig5} (see also
\cite{bl11} for more general discussions), or 
by the largest (in terms of corresponding eigenvalues)
 $L_i$ eigenvectors of $S^{+,\text{off}}\Psi_k^{\text{off}}= \left( \lambda_k^{\text{off}}\right)^{-1} A^{\text{off}} \Psi_k^{\text{off}}$.
Consequently, the use of Eq.~\eqref{offeig5} in constructing local
spaces will give a better approximation
compared to using multiple spectral problems that provides
a rate which is the product of $\prod_j 1/\Lambda_{*,j}$, where
$j$ represents the corresponding eigenvalue for $j$-th eigenvalue
problem. Consequently, if we set a threshold for the eigenvalue for 
each problem as $\Lambda_*$, then the convergence rate is $(1/\Lambda_*)^n$,
where $n$ is the number of eigenvalue problems are used. 
Note also that, each used eigenvalue problem will increase the dimension of 
the final reduced space.
 In general, one can show similar results for eigenvalue
problems considered earlier, as in Eqs.~\eqref{offeig2} and ~\eqref{offeig3}.

We present numerical results corresponding to the use 
of the eigenvalue problem in 
Eq.~\eqref{offeig5} in Table \ref{off20_parindep_Atp_S_offeig5}.
If we compare these results to Table \ref{off20_parindep_At_Sp}, we observe
that the convergence of the method is better than if the eigenvalue
problem of Eq.~\eqref{offeig3} is used.

\begin{table}[htb]
\centering

\begin{tabular}{|c|c|c|c|c|c|}
\hline 
\multirow{2}{*}{$\text{dim}(V_{\text{off}}^+)$} &\multirow{2}{*}{$\Lambda_*$} &
\multicolumn{2}{c|}{  $\|u-u^{+, \text{off}} \|$ (\%) }  \\
\cline{3-4} {}&{}&
$\hspace*{0.8cm}   L^{2}_\kappa(D)   \hspace*{0.8cm}$ &
$\hspace*{0.8cm}   H^{1}_\kappa(D)  \hspace*{0.8cm}$
\\
\hline\hline
       $1163$     &     $3.01{\times}10^{3}$   &  $3.11$    & $34.53$  \\
\hline
      $1524$    &   $5.35{\times}10^{3}$  & $0.25$    & $13.64$ \\
\hline
      $1885$    &   $7.51{\times}10^{3}$  & $0.20$    & $8.69$ \\
\hline
       $2607$   &  $7.76{\times}10^{3}$ & $0.12$  & $5.91$  \\

\hline
\end{tabular}

\caption{Relative errors between the fine scale solution and {\bf offline spaces} obtained from 
              using oversampled domains; Eigenvalue problem from
              Eq.~\eqref{offeig5}, $20{\times}20$ coarse mesh, harmonic snapshots, $\omega_i^+=3\times\omega_i$.}

\label{off20_parindep_Atp_S_offeig5}
\end{table}

\section{Conclusions}\label{sec:conclusions}

In this paper, we develop and investigate oversampling strategies for GMsFEM.
The GMsFEM offers a flexible framework for solving multiscale problems 
by constructing a reduced dimensional approximation for the solution
space.
In particular, GMsFEM constructs a local approximation space via appropriate
local spectral problems.
 We show that the use of oversampling strategies can yield a convergence
independent of the contrast and the small scales  under certain assumptions. 
The proof relies on
the fact that the local spectral problems that are used for basis construction
involve oversampled regions. The convergence of GMsFEM is proportional
to the maximum of the inverse of the eigenvalue such that the corresponding
eigenvector is not included in the coarse space. Our numerical results show
that the reciprocal of the eigenvalues decay faster when oversampling is used
(in particular, for the local spectral problem that is proposed in the paper).
We present some representative
numerical results where various
oversampling strategies are studied.  Our results compare the fine grid
solution with GMsFEM solution as well as the solution computed
in the snapshot space with GMsFEM solution.
We study the use of multiple spectral problems for enhanced accuracy
and discuss their relation to single spectral problems that use oversampled
regions where the latter provides an optimal space.
Both convergence analysis and numerical studies are presented.
Numerical results show that the proposed oversampling techniques
are efficient and have similar errors.
 We also present numerical
results for parameter-dependent problems using our proposed strategies.
The numerical results for each configuration are discussed in the paper.

\section*{Acknowledgments}

Y. Efendiev's work is
partially supported by the
DOE and NSF (DMS 0934837 and DMS 0811180).
J. Galvis would like to acknowledge partial support from DOE.
 This publication is based in part on work supported by Award
No. KUS-C1-016-04, made by King Abdullah University of Science
and Technology (KAUST).

\section*{References}

\bibliographystyle{plain}

\begin{thebibliography}{99}


\bibitem{aarnes04}
{ J. E.  Aarnes}, {\em On the use of a mixed multiscale finite element method
  for greater flexibility and increased speed or improved accuracy in reservoir
  simulation}, SIAM J. Multiscale Modeling and Simulation, 2 (2004),   421-439.


\bibitem{AKL}
{ J. E. Aarnes, S. Krogstad, and K.-A. Lie},
{\em A hierarchical multiscale method for two-phase flow based upon mixed finite elements and nonuniform grids},
  Multiscale Model. Simul. 5(2) (2006), pp. 337--363.


\bibitem{apwy07}
{  T. Arbogast, G. Pencheva, M. F. Wheeler, and I. Yotov}, {\em A multiscale mortar mixed finite element method}, SIAM J. Multiscale Modeling and Simulation, 6(1), 2007,   319-346.





\bibitem{bl11}
I.  Babu{\v{s}}ka and R. Lipton, {\em Optimal Local Approximation Spaces for Generalized Finite Element Methods with Application to Multiscale Problems}.  Multiscale Modeling and Simulation, SIAM 9 (2011) 373-406.





\bibitem{bo83}
I. Babu{\v{s}}ka and E. Osborn,
{\em Generalized Finite Element Methods: Their Performance and Their
Relation to Mixed Methods,}
SIAM J. Numer. Anal., 20 (1983), pp. 510-536.



\bibitem{bo10}
L. Berlyand and H. Owhadi, {\em Flux norm approach to finite dimensional homogenization approximations with non-separated scales and high contrast}, Arch. Ration. Mech. Anal., 198 (2010), pp. 677-721.



\bibitem{boyaval}
S. Boyaval, C. LeBris, T. Leli\`evre, Y. Maday, N. Nguyen, and A. Patera, {\em Reduced Basis Techniques for Stochastic Problems},
  Archives of Computational Methods in Engineering, 17:435-454, 2010.


\bibitem{cd07}
{ Y.~Chen and L.~Durlofsky}, {\em An ensemble level upscaling approach for efficient estimation of fine-scale production statistics using coarse-scale simulations}, SPE paper 106086, presented at the SPE Reservoir Simulation Symposium, Houston, Feb. 26-28 (2007).


\bibitem{cdgw03}
{  Y. Chen, L.J. Durlofsky, M. Gerritsen, and X.H. Wen}, {\em A coupled
  local-global upscaling approach for simulating flow in highly heterogeneous
  formations}, Advances in Water Resources, 26 (2003), pp. 1041--1060.




\bibitem{DGcoupling}
Y. Efendiev, J. Galvis, R. Lazarov, M. Moon, M. Sarkis,
{\it Generalized Multiscale Finite Element Method. Symmetric Interior Penalty Coupling}, 
Submitted.


\bibitem{eglw11}
Y. Efendiev, J. Galvis, R. Lazarov, and J. Willems, {\em Robust domain decomposition preconditioners for abstract symmetric positive definite bilinear forms}, ESAIM: Mathematical Modelling and Numerical Analysis, September 2012, 46, pp. 1175-1199. 


\bibitem{egh12}
{  Y. Efendiev, J. Galvis, and T. Hou}, {\em  Generalized Multiscale Finite Element Methods}, Submitted. http://arxiv.org/submit/631572.



\bibitem{egt11}
{  Y. Efendiev, J. Galvis, and F. Thomines}, {\em  A systematic coarse-scale model reduction technique
for parameter-dependent flows in highly heterogeneous media and its applications}, SIAM MMS 10(4), 1317-1343, 2012.


\bibitem{egw10}
{  Y. Efendiev, J. Galvis, and X. H. Wu}, {\em Multiscale finite element methods for high-contrast problems using local spectral basis functions}, Journal of Computational Physics. Volume 230, Issue 4, 20 February 2011, Pages 937-955.


\bibitem{eghe05}
{  Y. Efendiev, V. Ginting, T. Hou, and R. Ewing}, {\em Accurate multiscale finite element methods for two-phase flow simulations},
 J. Comp. Physics, 220 (1), pp. 155--174, 2006.


\bibitem{eh09}
{Y. Efendiev and T. Hou}, {\em Multiscale finite element methods. Theory
and applications}, Springer, 2009.


\bibitem{ehg04}
Y. Efendiev, T. Hou, and V. Ginting.
{\em Multiscale finite element methods for nonlinear problems and their
applications}. Comm. Math. Sci., 2:553�1�79, 2004.


\bibitem{ge09_3}
 J. Galvis and Y. Efendiev, {\it
 Domain decomposition  preconditioners for multiscale flows in high-contrast media: Reduced dimension coarse spaces}, SIAM J. Multiscale Modeling and Simulation, Volume 8, Issue 5,   1621-1644 (2010).





\bibitem{hw97}
{  T.Y. Hou and X.H. Wu}, {\em A multiscale finite element method for
  elliptic problems in composite materials and porous media}, Journal of
  Computational Physics, 134 (1997),   169-189.


\bibitem{hughes98}
{  T. Hughes, G. Feijoo, L. Mazzei, and J. Quincy}, {\em The variational
  multiscale method - a paradigm for computational mechanics}, Comput. Methods
  Appl. Mech. Engrg, 166 (1998),   3-24.


\bibitem{jennylt03}
{  P. Jenny, S.H. Lee, and H. Tchelepi}, {\em Multi-scale finite volume
  method for elliptic problems in subsurface flow simulation}, J. Comput.
  Phys., 187 (2003),   47-67.


\bibitem{GMNP07}
M.A. Grepl, Y. Maday, N.C. Nguyen, and A.T. Patera.
{\em Efficient reduced-basis treatment of non-affine and
 nonlinear partial differential equations}. ESIAM : M2AN, 41(2):575�1�75, 2007.




\bibitem{Patera}
G. Rozza, D. B. P Huynh, and A. T. Patera, {\em Reduced basis approximation and a posteriori error estimation for affinely parametrized elliptic coercive partial differential equations. Application to transport and continuum mechanics}. Arch Comput Methods Eng 15(3):229?275, 2008.



\end{thebibliography}

\def\cprime{$'$}

\end{document}